\documentclass[final,3p]{elsarticle}
\usepackage[latin1]{inputenc}
 \usepackage{graphics}
 \usepackage{graphicx}
 \usepackage{epsfig}
\usepackage{amssymb}
 \usepackage{amsthm}
 \usepackage{lineno}
 \usepackage{amsmath}
   \numberwithin{equation}{section}
\usepackage{mathrsfs}

\newtheorem{thm}{Theorem}[section]

\newtheorem{lem}[thm]{Lemma}

\newtheorem{defn}[thm]{Definition}

\biboptions{sort&compress,square}
\begin{document}
\begin{frontmatter}
\author{Tong Wu}
\ead{wut977@nenu.edu.cn}
\author{Yong Wang\corref{cor2}}
\ead{wangy581@nenu.edu.cn}
\cortext[cor2]{Corresponding author.}

\address{School of Mathematics and Statistics, Northeast Normal University,
Changchun, 130024, China}

\title{ Codazzi tensors and the quasi-statistical structure associated to affine connections on three-dimensional Lorentzian Lie groups}
\begin{abstract}
In this paper, we classify three-dimensional Lorentzian Lie groups on which Ricci tensors associated to Bott connections, canonical connections and Kobayashi-Nomizu connections are Codazzi tensors associated to these connections. We also classify three-dimensional Lorentzian Lie group with the quasi-statistical structure associated to Bott connections, canonical connections and Kobayashi-Nomizu connections.
\end{abstract}
\end{frontmatter}
\section{Introduction}
\indent In \cite{DC}, Andrzej and Shen studied  some geometric and topological consequences of the existence of a non-trivial Codazzi tensor on a Riemannian manifold. They also introduced Codazzi tensors associated to any linear connections. Bourguignon got the results of this tpye and gave the proof of the existence of such a tensor improses strong restrictions on the curvature operator in \cite{FB}. In \cite{MR}, Dajczer and Tojeiro found the correspondence between the Ribaucour transformation of a submanifold and Codazzi tensor exchanged with its second fundamental form. In \cite{JN}, authors defined a Codazzi tensor on conformally symmetric space, and characterized Einstein manifold and constant sectional curvature manifold by inequalities between certain functions about this tensor.\\
\indent In \cite{GM}, Merton and Gabe discussed the classification of Codazzi tensors with exactly two eigenfunctions on a Riemannian manifold of dimension three or higher.
   In \cite{A}, Blaga and Nannicini considered the statistical structure on a smooth manifold with a torsion-free affine connection, and they also gave the definition of the quasi-statistical structure, which is the generalization of the statistical structure. Wang gave algebraic Ricci solitons and affine Ricci solitons associated to canonical connections and Kobayashi-Nomizu connections on three-dimensional Lorentzian Lie groups respectively in \cite{wy,YY2}. In \cite{FB,J}, authors gave the definition of the Bott connection. In this paper, we classify three-dimensional Lorentzian Lie groups on which Ricci tensors associated to Bott connections, canonical connections and Kobayashi-Nomizu connections are Codazzi tensors associated to these connections. We also classify three-dimensional Lorentzian Lie group with the quasi-statistical structure associated to Bott connections, canonical connections and Kobayashi-Nomizu connections. \\
\indent In Section 2, we classify three-dimensional Lorentzian Lie groups on which Ricci tensors associated to Bott connections are Codazzi tensors associated to Bott connections. In Section 3, we classify three-dimensional Lorentzian Lie group with the quasi-statistical structure associated to Bott connections. In Section 4, we classify three-dimensional Lorentzian Lie groups on which Ricci tensors associated to canonical connections and Kobayashi-Nomizu connections are Codazzi tensors associated to canonical connections and Kobayashi-Nomizu connections. In Section 5, we classify three-dimensional Lorentzian Lie group with the quasi-statistical structure associated to canonical connections and Kobayashi-Nomizu connections.

\section{Codazzi tensors associated to Bott connections on three-dimensional Lorentzian Lie groups}
\indent Let $\{G_i\}_{i=1,\cdot\cdot\cdot,7}$, denote the connected, simply connected three-dimensional Lie group equipped with a left-invariant Lorentzian metric $g$ and having Lie algebra $\{\mathfrak{g}_i\}_{i=1,\cdot\cdot\cdot,7}$ and let $\nabla^L$ be the Levi-Civita connection of $G_i$. Nextly, we recall the definition of the Bott connection $\nabla^B$. Let $M$ be a smooth manifold, and let $TM=span\{\widetilde{e}_1,\widetilde{e}_2,\widetilde{e}_3\}$, then take the distribution $D=span\{\widetilde{e}_1,\widetilde{e}_2\}$ and $D^\bot=span\{\widetilde{e}_3\}$.\\
The definition of the Bott connection $\nabla^B$ is given as follows: (see \cite{FB}, \cite{J})
\begin{eqnarray}
\nabla^B_XY=
       \begin{cases}
        \pi_D(\nabla^L_XY),~~~&X,Y\in\Gamma^\infty(D) \\[2pt]
       \pi_D([X,Y]),~~~&X\in\Gamma^\infty(D^\bot),Y\in\Gamma^\infty(D)\\[2pt]
       \pi_{D^\bot}([X,Y]),~~~&X\in\Gamma^\infty(D),Y\in\Gamma^\infty(D^\bot)\\[2pt]
       \pi_{D^\bot}(\nabla^L_XY),~~~&X,Y\in\Gamma^\infty(D^\bot)\\[2pt]
       \end{cases}
\end{eqnarray}
where $\pi_D$(resp. $\pi_D^\bot$) the projection on $D$ (resp. $D^\bot$).\\
We define
\begin{equation}
R^B(X,Y)Z=\nabla^B_X\nabla^B_YZ-\nabla^B_Y\nabla^B_XZ-\nabla^B_{[X,Y]}Z.
\end{equation}
The Ricci tensor of $(G_i,g)$ associated to the Bott connection $\nabla^B$ is defined by
\begin{equation}
\rho^B(X,Y)=-g(R^B(X,\widetilde{e}_1)Y,\widetilde{e}_1)-g(R^B(X,\widetilde{e}_2)Y,\widetilde{e}_2)+g(R^B(X,\widetilde{e}_3)Y,\widetilde{e}_3),
\end{equation}
where ${\widetilde{e}_1,\widetilde{e}_2,\widetilde{e}_3}$ is a pseudo-orthonormal basis, with $\widetilde{e}_3$ timelike.\\
Let
\begin{equation}
\widetilde{\rho}^B(X,Y)=\frac{\rho^B(X,Y)+\rho^B(Y,X)}{2}.
\end{equation}
Let $\omega$ be a (0,2) tensor fileds, then
we define:
\begin{equation}
(\nabla_X\omega)(Y,Z):=X[\omega(Y,Z)]-\omega(\nabla_XY,Z)-\omega(Y,\nabla_XZ),
\end{equation}
for arbitrary vector fileds $X,Y,Z$.
\begin{defn} (\cite{DC}, $P_{17}$) Let $M$ be a smooth manifold endowed with a linear connection $\nabla$, the tensor fields $\omega$ is called a Codazzi tensor on $(M,\nabla)$, if it satisfies
\begin{equation}
f(X,Y,Z)=(\nabla_X\omega)(Y,Z)-(\nabla_Y\omega)(X,Z)=0,
\end{equation}
where $f$ is $C^\infty(M)$-linear for $X,Y,Z$, and $f(X,Y,Z)=-f(Y,X,Z)$.
\end{defn}
Then we have $\omega$ is a Codazzi tensor on $(G_i,\nabla)$ if and only if the following nine equations hold:
\begin{eqnarray}
       \begin{cases}
f(\widetilde{e}_1,\widetilde{e}_2,\widetilde{e}_j)=0\\[2pt]
f(\widetilde{e}_1,\widetilde{e}_3,\widetilde{e}_j)=0\\[2pt]
f(\widetilde{e}_2,\widetilde{e}_3,\widetilde{e}_j)=0\\[2pt]
 \end{cases}
\end{eqnarray}
where $1\leq j\leq3$.
\vskip 0.5 true cm
\noindent{\bf 2.1 Codazzi tensors of $G_1$}\\
\vskip 0.5 true cm
 By \cite{W}, we have the following Lie algebra of $G_1$ satisfies
\begin{equation}
[\widetilde{e}_1,\widetilde{e}_2]=\alpha\widetilde{e}_1-\beta\widetilde{e}_3,~~~[\widetilde{e}_1, \widetilde{e}_3]=-\alpha\widetilde{e}_1-\beta\widetilde{e}_2,~~~[\widetilde{e}_2, \widetilde{e}_3]=\beta\widetilde{e}_1+\alpha\widetilde{e}_2+\alpha\widetilde{e}_3,~\alpha\neq0.
\end{equation}
where ${\widetilde{e}_1,\widetilde{e}_2,\widetilde{e}_3}$ is a pseudo-orthonormal basis, with $\widetilde{e}_3$ timelike.
\begin{lem} The Bott connection $\nabla^B$ of $G_1$ is given by
\begin{align}
&\nabla^B_{\widetilde{e}_1}\widetilde{e}_1=-\alpha\widetilde{e}_2,~~~\nabla^B_{\widetilde{e}_1}\widetilde{e}_2=\alpha\widetilde{e}_1,~~~\nabla^B_{\widetilde{e}_1}\widetilde{e}_3=0,\nonumber\\
&\nabla^B_{\widetilde{e}_2}\widetilde{e}_1=0,~~~\nabla^B_{\widetilde{e}_2}\widetilde{e}_2=0,~~~\nabla^B_{\widetilde{e}_2}\widetilde{e}_3=\alpha\widetilde{e}_3,\nonumber\\
&\nabla^B_{\widetilde{e}_3}\widetilde{e}_1=\alpha\widetilde{e}_1+\beta\widetilde{e}_2,~~~\nabla^B_{\widetilde{e}_3}\widetilde{e}_2=-\beta\widetilde{e}_1-\alpha\widetilde{e}_2,~~~\nabla^B_{\widetilde{e}_3}\widetilde{e}_3=0.
\end{align}
\end{lem}
\begin{lem} The curvature $R^B$ of the Bott connection $\nabla^B$ of $(G_1,g)$ is given by
\begin{align}
&R^B(\widetilde{e}_1,\widetilde{e}_2)\widetilde{e}_1=\alpha\beta\widetilde{e}_1+(\alpha^2+\beta^2)\widetilde{e}_2,~~~R^B(\widetilde{e}_1,\widetilde{e}_2)\widetilde{e}_2=-(\alpha^2+\beta^2)\widetilde{e}_1-\alpha\beta\widetilde{e}_2,~~~R^B(\widetilde{e}_1,\widetilde{e}_2)\widetilde{e}_3=0,\nonumber\\
&R^B(\widetilde{e}_1,\widetilde{e}_3)\widetilde{e}_1=-3\alpha^2\widetilde{e}_2,~~~R^B(\widetilde{e}_1,\widetilde{e}_3)\widetilde{e}_2=-\alpha^2\widetilde{e}_1,~~~R^B(\widetilde{e}_1,\widetilde{e}_3)\widetilde{e}_3=\alpha\beta\widetilde{e}_3,\nonumber\\
&R^B(\widetilde{e}_2,\widetilde{e}_3)\widetilde{e}_1=-\alpha^2\widetilde{e}_1,~~~R^B(\widetilde{e}_2,\widetilde{e}_3)\widetilde{e}_2=\alpha^2\widetilde{e}_2,~~~R^B(\widetilde{e}_2,\widetilde{e}_3)\widetilde{e}_3=-\alpha^2\widetilde{e}_3.
\end{align}
\end{lem}
By (2.3), we have
\begin{align}
&\rho^B(\widetilde{e}_1,\widetilde{e}_1)=-(\alpha^2+\beta^2),~~~\rho^B(\widetilde{e}_1,\widetilde{e}_2)=\alpha\beta,~~~\rho^B(\widetilde{e}_1,\widetilde{e}_3)=-\alpha\beta,\nonumber\\
&\rho^B(\widetilde{e}_2,\widetilde{e}_1)=\alpha\beta,~~~\rho^B(\widetilde{e}_2,\widetilde{e}_2)=-(\alpha^2+\beta^2),~~~\rho^B(\widetilde{e}_2,\widetilde{e}_3)=\alpha^2,\nonumber\\
&\rho^B(\widetilde{e}_3,\widetilde{e}_1)=\rho^B(\widetilde{e}_3,\widetilde{e}_2)=\rho^B(\widetilde{e}_3,\widetilde{e}_3)=0.
\end{align}
Then,
\begin{align}
&\widetilde{\rho}^B(\widetilde{e}_1,\widetilde{e}_1)=-(\alpha^2+\beta^2),~~~\widetilde{\rho}^B(\widetilde{e}_1,\widetilde{e}_2)=\alpha\beta,~~~\widetilde{\rho}^B(\widetilde{e}_1,\widetilde{e}_3)=-\frac{\alpha\beta}{2},\nonumber\\
&\widetilde{\rho}^B(\widetilde{e}_2,\widetilde{e}_2)=-(\alpha^2+\beta^2),~~~\widetilde{\rho}^B(\widetilde{e}_2,\widetilde{e}_3)=\frac{\alpha^2}{2},~~~\widetilde{\rho}^B(\widetilde{e}_3,\widetilde{e}_3)=0.
\end{align}
By (2.5), we have
\begin{align}
&(\nabla^B_{\widetilde{e}_1}\widetilde{\rho}^B)(\widetilde{e}_2,\widetilde{e}_1)=0,~~~(\nabla^B_{\widetilde{e}_2}\widetilde{\rho}^B)(\widetilde{e}_1,\widetilde{e}_1)=0,~~~(\nabla^B_{\widetilde{e}_1}\widetilde{\rho}^B)(\widetilde{e}_2,\widetilde{e}_2)=-2\alpha^2\beta,\nonumber\\
&(\nabla^B_{\widetilde{e}_2}\widetilde{\rho}^B)(\widetilde{e}_1,\widetilde{e}_2)=0,~~~(\nabla^B_{\widetilde{e}_1}\widetilde{\rho}^B)(\widetilde{e}_2,\widetilde{e}_3)=\frac{\alpha^2\beta}{2},~~~(\nabla^B_{\widetilde{e}_2}\widetilde{\rho}^B)(\widetilde{e}_1,\widetilde{e}_3)=\frac{\alpha^2\beta}{2},\nonumber\\
&(\nabla^B_{\widetilde{e}_1}\widetilde{\rho}^B)(\widetilde{e}_3,\widetilde{e}_1)=\frac{\alpha^3}{2},~~~(\nabla^B_{\widetilde{e}_3}\widetilde{\rho}^B)(\widetilde{e}_1,\widetilde{e}_1)=2\alpha^3,~~~(\nabla^B_{\widetilde{e}_1}\widetilde{\rho}^B)(\widetilde{e}_3,\widetilde{e}_2)=\frac{\alpha^2\beta}{2},\nonumber\\
&(\nabla^B_{\widetilde{e}_3}\widetilde{\rho}^B)(\widetilde{e}_1,\widetilde{e}_2)=0,~~~(\nabla^B_{\widetilde{e}_1}\widetilde{\rho}^B)(\widetilde{e}_3,\widetilde{e}_3)=0,~~~(\nabla^B_{\widetilde{e}_3}\widetilde{\rho}^B)(\widetilde{e}_1,\widetilde{e}_3)=0,\nonumber\\
&(\nabla^B_{\widetilde{e}_2}\widetilde{\rho}^B)(\widetilde{e}_3,\widetilde{e}_1)=\frac{\alpha^2\beta}{2},~~~(\nabla^B_{\widetilde{e}_3}\widetilde{\rho}^B)(\widetilde{e}_2,\widetilde{e}_1)=0,~~~(\nabla^B_{\widetilde{e}_2}\widetilde{\rho}^B)(\widetilde{e}_3,\widetilde{e}_2)=-\frac{\alpha^2}{2},\nonumber\\
&(\nabla^B_{\widetilde{e}_3}\widetilde{\rho}^B)(\widetilde{e}_2,\widetilde{e}_2)=-2\alpha^2,~~~(\nabla^B_{\widetilde{e}_2}\widetilde{\rho}^B)(\widetilde{e}_3,\widetilde{e}_3)=0,~~~(\nabla^B_{\widetilde{e}_3}\widetilde{\rho}^B)(\widetilde{e}_2,\widetilde{e}_3)=\frac{\alpha}{2}(\alpha^2-\beta^2).
\end{align}
Then, if $\widetilde{\rho}^B$  is a Codazzi tensor on $(G_1,\nabla^B)$, by (2.6) and (2.7), we have the following three equations:
\begin{eqnarray}
       \begin{cases}
        2\alpha^2\beta=0 \\[2pt]
      \frac{3\alpha^3}{2}=0\\[2pt]
       \frac{\alpha}{2}(\alpha^2-\beta^2)=0\\[2pt]
       \end{cases}
\end{eqnarray}
By solving (2.14) , we get $\alpha=0$, there is a contradiction. So\\
\begin{thm}
$\widetilde{\rho}^B$  is not a Codazzi tensor on $(G_1,\nabla^B)$.
\end{thm}
\vskip 0.5 true cm
\noindent{\bf 2.2 Codazzi tensors of $G_2$}\\
\vskip 0.5 true cm
\indent By \cite{W}, we have the following Lie algebra of $G_2$ satisfies
\begin{equation}
[\widetilde{e}_1, \widetilde{e}_2]=\gamma\widetilde{e}_2-\beta\widetilde{e}_3,~~~[\widetilde{e}_1, \widetilde{e}_3]=-\beta\widetilde{e}_2-\gamma\widetilde{e}_3,~~~[\widetilde{e}_2, \widetilde{e}_3]=\alpha\widetilde{e}_1,~\gamma\neq0.
\end{equation}
where ${\widetilde{e}_1,\widetilde{e}_2,\widetilde{e}_3}$ is a pseudo-orthonormal basis, with $\widetilde{e}_3$ timelike.
\begin{lem} The Bott connection $\nabla^B$ of $G_2$ is given by
\begin{align}
&\nabla^B_{\widetilde{e}_1}\widetilde{e}_1=0,~~~\nabla^B_{\widetilde{e}_1}\widetilde{e}_2=0,~~~\nabla^B_{\widetilde{e}_1}\widetilde{e}_3=-\gamma\widetilde{e}_3,\nonumber\\
&\nabla^B_{\widetilde{e}_2}\widetilde{e}_1=-\gamma\widetilde{e}_2,~~~\nabla^B_{\widetilde{e}_2}\widetilde{e}_2=\gamma\widetilde{e}_1,~~~\nabla^B_{\widetilde{e}_2}\widetilde{e}_3=0,\nonumber\\
&\nabla^B_{\widetilde{e}_3}\widetilde{e}_1=\beta\widetilde{e}_2,~~~\nabla^B_{\widetilde{e}_3}\widetilde{e}_2=-\alpha\widetilde{e}_1,~~~\nabla^B_{\widetilde{e}_3}\widetilde{e}_3=0.
\end{align}
\end{lem}
\begin{lem} The curvature $R^B$ of the Bott connection $\nabla^B$ of $(G_2,g)$ is given by
\begin{align}
&R^B(\widetilde{e}_1,\widetilde{e}_2)\widetilde{e}_1=(\beta^2+\gamma^2)\widetilde{e}_2,~~~R^B(\widetilde{e}_1,\widetilde{e}_2)\widetilde{e}_2=-(\gamma^2+\alpha\beta)\widetilde{e}_1,~~~R^B(\widetilde{e}_1,\widetilde{e}_2)\widetilde{e}_3=0,\nonumber\\
&R^B(\widetilde{e}_1,\widetilde{e}_3)\widetilde{e}_1=0,~~~R^B(\widetilde{e}_1,\widetilde{e}_3)\widetilde{e}_2=\gamma(\alpha-\beta)\widetilde{e}_1,~~~R^B(\widetilde{e}_1,\widetilde{e}_3)\widetilde{e}_3=0,\nonumber\\
&R^B(\widetilde{e}_2,\widetilde{e}_3)\widetilde{e}_1=\gamma(\beta-\alpha)\widetilde{e}_1,~~~R^B(\widetilde{e}_2,\widetilde{e}_3)\widetilde{e}_2=\gamma(\alpha-\beta)\widetilde{e}_2,~~~R^B(\widetilde{e}_2,\widetilde{e}_3)\widetilde{e}_3=\alpha\gamma\widetilde{e}_3.
\end{align}
\end{lem}
By (2.3), we have
\begin{align}
&\rho^B(\widetilde{e}_1,\widetilde{e}_1)=-(\beta^2+\gamma^2),~~~\rho^B(\widetilde{e}_1,\widetilde{e}_2)=0,~~~\rho^B(\widetilde{e}_1,\widetilde{e}_3)=0,\nonumber\\
&\rho^B(\widetilde{e}_2,\widetilde{e}_1)=0,~~~\rho^B(\widetilde{e}_2,\widetilde{e}_2)=-(\gamma^2+\alpha\beta),~~~\rho^B(\widetilde{e}_2,\widetilde{e}_3)=-\alpha\gamma,\nonumber\\
&\rho^B(\widetilde{e}_3,\widetilde{e}_1)=\rho^B(\widetilde{e}_3,\widetilde{e}_2)=\rho^B(\widetilde{e}_3,\widetilde{e}_3)=0.
\end{align}
Then,
\begin{align}
&\widetilde{\rho}^B(\widetilde{e}_1,\widetilde{e}_1)=-(\beta^2+\gamma^2),~~~\widetilde{\rho}^B(\widetilde{e}_1,\widetilde{e}_2)=0,~~~\widetilde{\rho}^B(\widetilde{e}_1,\widetilde{e}_3)=0,\nonumber\\
&\widetilde{\rho}^B(\widetilde{e}_2,\widetilde{e}_2)=-(\gamma^2+\alpha\beta),~~~\widetilde{\rho}^B(\widetilde{e}_2,\widetilde{e}_3)=-\frac{\alpha\gamma}{2},~~~\widetilde{\rho}^B(\widetilde{e}_3,\widetilde{e}_3)=0.
\end{align}
By (2.5), we have
\begin{align}
&(\nabla^B_{\widetilde{e}_1}\widetilde{\rho}^B)(\widetilde{e}_2,\widetilde{e}_1)=0,~~~(\nabla^B_{\widetilde{e}_2}\widetilde{\rho}^B)(\widetilde{e}_1,\widetilde{e}_1)=0,~~~(\nabla^B_{\widetilde{e}_1}\widetilde{\rho}^B)(\widetilde{e}_2,\widetilde{e}_2)=0,\nonumber\\
&(\nabla^B_{\widetilde{e}_2}\widetilde{\rho}^B)(\widetilde{e}_1,\widetilde{e}_2)=\gamma(\beta^2-\alpha\beta),~~~(\nabla^B_{\widetilde{e}_1}\widetilde{\rho}^B)(\widetilde{e}_2,\widetilde{e}_3)=-\frac{\alpha\gamma^2}{2},~~~(\nabla^B_{\widetilde{e}_2}\widetilde{\rho}^B)(\widetilde{e}_1,\widetilde{e}_3)=-\frac{\alpha\gamma^2}{2},\nonumber\\
&(\nabla^B_{\widetilde{e}_1}\widetilde{\rho}^B)(\widetilde{e}_3,\widetilde{e}_1)=0,~~~(\nabla^B_{\widetilde{e}_3}\widetilde{\rho}^B)(\widetilde{e}_1,\widetilde{e}_1)=0,~~~(\nabla^B_{\widetilde{e}_1}\widetilde{\rho}^B)(\widetilde{e}_3,\widetilde{e}_2)=-\frac{\alpha\gamma^2}{2},\nonumber\\
&(\nabla^B_{\widetilde{e}_3}\widetilde{\rho}^B)(\widetilde{e}_1,\widetilde{e}_2)=\gamma^2(\beta-\alpha),~~~(\nabla^B_{\widetilde{e}_1}\widetilde{\rho}^B)(\widetilde{e}_3,\widetilde{e}_3)=0,~~~(\nabla^B_{\widetilde{e}_3}\widetilde{\rho}^B)(\widetilde{e}_1,\widetilde{e}_3)=\frac{\alpha\beta\gamma}{2},\nonumber\\
&(\nabla^B_{\widetilde{e}_2}\widetilde{\rho}^B)(\widetilde{e}_3,\widetilde{e}_1)=-\frac{\alpha\gamma^2}{2},~~~(\nabla^B_{\widetilde{e}_3}\widetilde{\rho}^B)(\widetilde{e}_2,\widetilde{e}_1)=\gamma^2(\beta-\alpha),~~~(\nabla^B_{\widetilde{e}_2}\widetilde{\rho}^B)(\widetilde{e}_3,\widetilde{e}_2)=0,\nonumber\\
&(\nabla^B_{\widetilde{e}_3}\widetilde{\rho}^B)(\widetilde{e}_2,\widetilde{e}_2)=0,~~~(\nabla^B_{\widetilde{e}_2}\widetilde{\rho}^B)(\widetilde{e}_3,\widetilde{e}_3)=0,~~~(\nabla^B_{\widetilde{e}_3}\widetilde{\rho}^B)(\widetilde{e}_2,\widetilde{e}_3)=0.
\end{align}
Then, if $\widetilde{\rho}^B$  is a Codazzi tensor on $(G_2,\nabla^B)$, by (2.6) and (2.7), we have the following three equations:
\begin{eqnarray}
       \begin{cases}
        \gamma(\beta^2-\alpha\beta)=0 \\[2pt]
      \frac{\alpha\beta\gamma}{2}=0\\[2pt]
       \gamma^2(\frac{\alpha}{2}-\beta)=0\\[2pt]
       \end{cases}
\end{eqnarray}
By solving (2.21), we get $\alpha=\beta=0$, there is a contradiction. So\\
\begin{thm}
 $\widetilde{\rho}^B$  is not a Codazzi tensor on $(G_2,\nabla^B)$.
\end{thm}
\vskip 0.5 true cm
\noindent{\bf 2.3 Codazzi tensors of $G_3$}\\
\vskip 0.5 true cm
\indent By \cite{W}, we have the following Lie algebra of $G_3$ satisfies
\begin{equation}
[\widetilde{e}_1, \widetilde{e}_2]=-\gamma\widetilde{e}_3,~~~[\widetilde{e}_1, \widetilde{e}_3]=-\beta\widetilde{e}_2,~~~[\widetilde{e}_2, \widetilde{e}_3]=\alpha\widetilde{e}_1.
\end{equation}
where ${\widetilde{e}_1,\widetilde{e}_2,\widetilde{e}_3}$ is a pseudo-orthonormal basis, with $\widetilde{e}_3$ timelike.
\begin{lem} The Bott connection $\nabla^B$ of $G_3$ is given by
\begin{align}
&\nabla^B_{\widetilde{e}_1}\widetilde{e}_1=0,~~~\nabla^B_{\widetilde{e}_1}\widetilde{e}_2=0,~~~\nabla^B_{\widetilde{e}_1}\widetilde{e}_3=-\gamma\widetilde{e}_3,\nonumber\\
&\nabla^B_{\widetilde{e}_2}\widetilde{e}_1=0,~~~\nabla^B_{\widetilde{e}_2}\widetilde{e}_2=0,~~~\nabla^B_{\widetilde{e}_2}\widetilde{e}_3=0,\nonumber\\
&\nabla^B_{\widetilde{e}_3}\widetilde{e}_1=\beta\widetilde{e}_2,~~~\nabla^B_{\widetilde{e}_3}\widetilde{e}_2=-\alpha\widetilde{e}_1,~~~\nabla^B_{\widetilde{e}_3}\widetilde{e}_3=0.
\end{align}
\end{lem}
\begin{lem} The curvature $R^B$ of the Bott connection $\nabla^B$ of $(G_3,g)$ is given by
\begin{align}
&R^B(\widetilde{e}_1,\widetilde{e}_2)\widetilde{e}_1=\beta\gamma\widetilde{e}_2,~~~R^B(\widetilde{e}_1,\widetilde{e}_2)\widetilde{e}_2=-\alpha\gamma\widetilde{e}_1,~~~R^B(\widetilde{e}_1,\widetilde{e}_2)\widetilde{e}_3=0,\nonumber\\
&R^B(\widetilde{e}_1,\widetilde{e}_3)\widetilde{e}_1=0,~~~R^B(\widetilde{e}_1,\widetilde{e}_3)\widetilde{e}_2=0,~~~R^B(\widetilde{e}_1,\widetilde{e}_3)\widetilde{e}_3=0,\nonumber\\
&R^B(\widetilde{e}_2,\widetilde{e}_3)\widetilde{e}_1=0,~~~R^B(\widetilde{e}_2,\widetilde{e}_3)\widetilde{e}_2=0,~~~R^B(\widetilde{e}_2,\widetilde{e}_3)\widetilde{e}_3=0.
\end{align}
\end{lem}
By (2.3), we have
\begin{align}
&\rho^B(\widetilde{e}_1,\widetilde{e}_1)=-\beta\gamma,~~~\rho^B(\widetilde{e}_1,\widetilde{e}_2)=0,~~~\rho^B(\widetilde{e}_1,\widetilde{e}_3)=0,\nonumber\\
&\rho^B(\widetilde{e}_2,\widetilde{e}_1)=0,~~~\rho^B(\widetilde{e}_2,\widetilde{e}_2)=-\alpha\gamma,~~~\rho^B(\widetilde{e}_2,\widetilde{e}_3)=0,\nonumber\\
&\rho^B(\widetilde{e}_3,\widetilde{e}_1)=\rho^B(\widetilde{e}_3,\widetilde{e}_2)=\rho^B(\widetilde{e}_3,\widetilde{e}_3)=0.
\end{align}
Then,
\begin{align}
&\widetilde{\rho}^B(\widetilde{e}_1,\widetilde{e}_1)=-\beta\gamma,~~~\widetilde{\rho}^B(\widetilde{e}_1,\widetilde{e}_2)=\widetilde{\rho}^B(\widetilde{e}_1,\widetilde{e}_3)=0,\nonumber\\
&\widetilde{\rho}^B(\widetilde{e}_2,\widetilde{e}_2)=-\alpha\gamma,~~~\widetilde{\rho}^B(\widetilde{e}_2,\widetilde{e}_3)=\widetilde{\rho}^B(\widetilde{e}_3,\widetilde{e}_3)=0.
\end{align}
By (2.5), we have
\begin{align}
&(\nabla^B_{\widetilde{e}_1}\widetilde{\rho}^B)(\widetilde{e}_2,\widetilde{e}_j)=(\nabla^B_{\widetilde{e}_2}\widetilde{\rho}^B)(\widetilde{e}_1,\widetilde{e}_j)=(\nabla^B_{\widetilde{e}_1}\widetilde{\rho}^B)(\widetilde{e}_3,\widetilde{e}_j)=0,\nonumber\\
&(\nabla^B_{\widetilde{e}_3}\widetilde{\rho}^B)(\widetilde{e}_1,\widetilde{e}_j)=(\nabla^B_{\widetilde{e}_2}\widetilde{\rho}^B)(\widetilde{e}_3,\widetilde{e}_j)=(\nabla^B_{\widetilde{e}_3}\widetilde{\rho}^B)(\widetilde{e}_2,\widetilde{e}_j)=0,
\end{align}
where $1\leq j\leq3$.\\
Then, we get
\begin{thm}
$\widetilde{\rho}^B$ is a Codazzi tensor on $(G_3,\nabla^B)$.
\end{thm}
\vskip 0.5 true cm
\noindent{\bf 2.4 Codazzi tensors of $G_4$}\\
\vskip 0.5 true cm
\indent By \cite{W}, we have the following Lie algebra of $G_4$ satisfies
\begin{equation}
[\widetilde{e}_1, \widetilde{e}_2]=-\widetilde{e}_2+(2\eta-\beta)\widetilde{e}_3,~~\eta=\pm 1,~~~[\widetilde{e}_1, \widetilde{e}_3]=-\beta\widetilde{e}_2+\widetilde{e}_3,~~~[\widetilde{e}_2, \widetilde{e}_3]=\alpha\widetilde{e}_1.
\end{equation}
where ${\widetilde{e}_1,\widetilde{e}_2,\widetilde{e}_3}$ is a pseudo-orthonormal basis, with $\widetilde{e}_3$ timelike.
\begin{lem} The Bott connection $\nabla^B$ of $G_4$ is given by
\begin{align}
&\nabla^B_{\widetilde{e}_1}\widetilde{e}_1=0,~~~\nabla^B_{\widetilde{e}_1}\widetilde{e}_2=0,~~~\nabla^B_{\widetilde{e}_1}\widetilde{e}_3=\widetilde{e}_3,\nonumber\\
&\nabla^B_{\widetilde{e}_2}\widetilde{e}_1=\widetilde{e}_2,~~~\nabla^B_{\widetilde{e}_2}\widetilde{e}_2=-\widetilde{e}_1,~~~\nabla^B_{\widetilde{e}_2}\widetilde{e}_3=0,\nonumber\\
&\nabla^B_{\widetilde{e}_3}\widetilde{e}_1=\beta\widetilde{e}_2,~~~\nabla^B_{\widetilde{e}_3}\widetilde{e}_2=-\alpha\widetilde{e}_1,~~~\nabla^B_{\widetilde{e}_3}\widetilde{e}_3=0.
\end{align}
\end{lem}
\begin{lem} The curvature $R^B$ of the Bott connection $\nabla^B$ of $(G_4,g)$ is given by
\begin{align}
&R^B(\widetilde{e}_1,\widetilde{e}_2)\widetilde{e}_1=(\beta-\eta)^2\widetilde{e}_2,~~~R^B(\widetilde{e}_1,\widetilde{e}_2)\widetilde{e}_2=(2\alpha\eta-\alpha\beta-1)\widetilde{e}_1,~~~R^B(\widetilde{e}_1,\widetilde{e}_2)\widetilde{e}_3=0,\nonumber\\
&R^B(\widetilde{e}_1,\widetilde{e}_3)\widetilde{e}_1=0,~~~R^B(\widetilde{e}_1,\widetilde{e}_3)\widetilde{e}_2=(\alpha-\beta)\widetilde{e}_1,~~~R^B(\widetilde{e}_1,\widetilde{e}_3)\widetilde{e}_3=0,\nonumber\\
&R^B(\widetilde{e}_2,\widetilde{e}_3)\widetilde{e}_1=(\alpha-\beta)\widetilde{e}_1,~~~R^B(\widetilde{e}_2,\widetilde{e}_3)\widetilde{e}_2=(\beta-\alpha)\widetilde{e}_2,~~~R^B(\widetilde{e}_2,\widetilde{e}_3)\widetilde{e}_3=-\alpha\widetilde{e}_3.
\end{align}
\end{lem}
By (2.3), we have
\begin{align}
&\rho^B(\widetilde{e}_1,\widetilde{e}_1)=-(\beta-\eta)^2,~~~\rho^B(\widetilde{e}_1,\widetilde{e}_2)=0,~~~\rho^B(\widetilde{e}_1,\widetilde{e}_3)=0,\nonumber\\
&\rho^B(\widetilde{e}_2,\widetilde{e}_1)=(2\alpha\eta-\alpha\beta-1),~~~\rho^B(\widetilde{e}_2,\widetilde{e}_2)=\alpha,~~~\rho^B(\widetilde{e}_2,\widetilde{e}_3)=0,\nonumber\\
&\rho^B(\widetilde{e}_3,\widetilde{e}_1)=\rho^B(\widetilde{e}_3,\widetilde{e}_2)=\rho^B(\widetilde{e}_3,\widetilde{e}_3)=0.
\end{align}
Then,
\begin{align}
&\widetilde{\rho}^B(\widetilde{e}_1,\widetilde{e}_1)=-(\beta-\eta)^2,~~~\widetilde{\rho}^B(\widetilde{e}_1,\widetilde{e}_2)=0,~~~\widetilde{\rho}^B(\widetilde{e}_1,\widetilde{e}_3)=0,\nonumber\\
&\widetilde{\rho}^B(\widetilde{e}_2,\widetilde{e}_2)=(2\alpha\eta-\alpha\beta-1),~~~\widetilde{\rho}^B(\widetilde{e}_2,\widetilde{e}_3)=\frac{\alpha}{2},~~~\widetilde{\rho}^B(\widetilde{e}_3,\widetilde{e}_3)=0.
\end{align}
By (2.5), we have
\begin{align}
&(\nabla^B_{\widetilde{e}_1}\widetilde{\rho}^B)(\widetilde{e}_2,\widetilde{e}_1)=0,~~~(\nabla^B_{\widetilde{e}_2}\widetilde{\rho}^B)(\widetilde{e}_1,\widetilde{e}_1)=0,~~~(\nabla^B_{\widetilde{e}_1}\widetilde{\rho}^B)(\widetilde{e}_2,\widetilde{e}_2)=\alpha\beta+2\beta\eta-2\alpha\eta-\beta^2,\nonumber\\
&(\nabla^B_{\widetilde{e}_2}\widetilde{\rho}^B)(\widetilde{e}_1,\widetilde{e}_2)=-\frac{\alpha}{2},~~~(\nabla^B_{\widetilde{e}_1}\widetilde{\rho}^B)(\widetilde{e}_2,\widetilde{e}_3)=-\frac{\alpha}{2},~~~(\nabla^B_{\widetilde{e}_2}\widetilde{\rho}^B)(\widetilde{e}_1,\widetilde{e}_3)=0,\nonumber\\
&(\nabla^B_{\widetilde{e}_1}\widetilde{\rho}^B)(\widetilde{e}_3,\widetilde{e}_1)=0,~~~(\nabla^B_{\widetilde{e}_3}\widetilde{\rho}^B)(\widetilde{e}_1,\widetilde{e}_1)=0,~~~(\nabla^B_{\widetilde{e}_1}\widetilde{\rho}^B)(\widetilde{e}_3,\widetilde{e}_2)=-\frac{\alpha}{2},\nonumber\\
&(\nabla^B_{\widetilde{e}_3}\widetilde{\rho}^B)(\widetilde{e}_1,\widetilde{e}_2)=\beta-\alpha,~~~(\nabla^B_{\widetilde{e}_1}\widetilde{\rho}^B)(\widetilde{e}_3,\widetilde{e}_3)=0,~~~(\nabla^B_{\widetilde{e}_3}\widetilde{\rho}^B)(\widetilde{e}_1,\widetilde{e}_3)=-\frac{\alpha\beta}{2},\nonumber\\
&(\nabla^B_{\widetilde{e}_2}\widetilde{\rho}^B)(\widetilde{e}_3,\widetilde{e}_1)=-\frac{\alpha}{2},~~~(\nabla^B_{\widetilde{e}_3}\widetilde{\rho}^B)(\widetilde{e}_2,\widetilde{e}_1)=\beta-\alpha,~~~(\nabla^B_{\widetilde{e}_2}\widetilde{\rho}^B)(\widetilde{e}_3,\widetilde{e}_2)=0,\nonumber\\
&(\nabla^B_{\widetilde{e}_3}\widetilde{\rho}^B)(\widetilde{e}_2,\widetilde{e}_2)=0,~~~(\nabla^B_{\widetilde{e}_2}\widetilde{\rho}^B)(\widetilde{e}_3,\widetilde{e}_3)=0,~~~(\nabla^B_{\widetilde{e}_3}\widetilde{\rho}^B)(\widetilde{e}_2,\widetilde{e}_3)=0.
\end{align}
Then, if $\widetilde{\rho}^B$ is a Codazzi tensor on $(G_4,\nabla^B)$ , by (2.6) and (2.7), we have the following three equations:
\begin{eqnarray}
       \begin{cases}
        (\beta-\eta)^2+2\alpha\eta-\alpha\beta-1=0 \\[2pt]
      \frac{\alpha}{2}-\beta=0\\[2pt]
       \frac{\alpha\beta}{2}=0\\[2pt]
       \end{cases}
\end{eqnarray}
By solving (2.34), we get
\begin{thm}
$\widetilde{\rho}^B$ is a Codazzi tensor on $(G_4,\nabla^B)$ if and only if $\alpha=\beta=0$.
\end{thm}
\noindent{\bf 2.5 Codazzi tensors of $G_5$}\\
\vskip 0.5 true cm
\indent By \cite{W}, we have the following Lie algebra of $G_5$ satisfies
\begin{equation}
[\widetilde{e}_1, \widetilde{e}_2]=0,~~~[\widetilde{e}_1, \widetilde{e}_3]=\alpha\widetilde{e}_1+\beta\widetilde{e}_2,~~~[\widetilde{e}_2, \widetilde{e}_3]=\gamma\widetilde{e}_1+\delta\widetilde{e}_2,~~~\alpha+\delta\neq 0,~~~\alpha\gamma+\beta\delta=0.
\end{equation}
where ${\widetilde{e}_1,\widetilde{e}_2,\widetilde{e}_3}$ is a pseudo-orthonormal basis, with $\widetilde{e}_3$ timelike.
\begin{lem} The Bott connection $\nabla^B$ of $G_5$ is given by
\begin{align}
&\nabla^B_{\widetilde{e}_1}\widetilde{e}_1=0,~~~\nabla^B_{\widetilde{e}_1}\widetilde{e}_2=0,~~~\nabla^B_{\widetilde{e}_1}\widetilde{e}_3=0,\nonumber\\
&\nabla^B_{\widetilde{e}_2}\widetilde{e}_1=0,~~~\nabla^B_{\widetilde{e}_2}\widetilde{e}_2=0,~~~\nabla^B_{\widetilde{e}_2}\widetilde{e}_3=0,\nonumber\\
&\nabla^B_{\widetilde{e}_3}\widetilde{e}_1=-\alpha\widetilde{e}_1-\beta\widetilde{e}_2,~~~\nabla^B_{\widetilde{e}_3}\widetilde{e}_2=-\gamma\widetilde{e}_1-\delta\widetilde{e}_2,~~~\nabla^B_{\widetilde{e}_3}\widetilde{e}_3=0.
\end{align}
\end{lem}
\begin{lem} The curvature $R^B$ of the Bott connection $\nabla^B$ of $(G_5,g)$ is given by
\begin{align}
R^B(\widetilde{e}_i,\widetilde{e}_j)\widetilde{e}_k=0,
\end{align}
for any $(i,j,k)$.
\end{lem}
By (2.3), we have
\begin{align}
\rho^B(\widetilde{e}_i,\widetilde{e}_j)=0,
\end{align}
then,
\begin{align}
&\widetilde{\rho}^B(\widetilde{e}_i,\widetilde{e}_j)=0,
\end{align}
for any pairs $(i,j)$.\\
By (2.5), we have\\
\begin{align}
&(\nabla^B_{\widetilde{e}_1}\widetilde{\rho}^B)(\widetilde{e}_2,\widetilde{e}_j)=(\nabla^B_{\widetilde{e}_2}\widetilde{\rho}^B)(\widetilde{e}_1,\widetilde{e}_j)=(\nabla^B_{\widetilde{e}_1}\widetilde{\rho}^B)(\widetilde{e}_3,\widetilde{e}_j)=0,\nonumber\\
&(\nabla^B_{\widetilde{e}_3}\widetilde{\rho}^B)(\widetilde{e}_1,\widetilde{e}_j)=(\nabla^B_{\widetilde{e}_2}\widetilde{\rho}^B)(\widetilde{e}_3,\widetilde{e}_j)=(\nabla^B_{\widetilde{e}_3}\widetilde{\rho}^B)(\widetilde{e}_2,\widetilde{e}_j)=0,
\end{align}
where $1\leq j\leq3$.\\
Then, we get
\begin{thm}
$\widetilde{\rho}^B$ is a Codazzi tensor on $(G_5,\nabla^B)$.
\end{thm}
\vskip 0.5 true cm
\noindent{\bf 2.6 Codazzi tensors of $G_6$}\\
\vskip 0.5 true cm
\indent By \cite{W}, we have the following Lie algebra of $G_6$ satisfies
\begin{equation}
[\widetilde{e}_1, \widetilde{e}_2]=\alpha\widetilde{e}_2+\beta\widetilde{e}_3,~~~[\widetilde{e}_1, \widetilde{e}_3]=\gamma\widetilde{e}_2+\delta\widetilde{e}_3,~~~[\widetilde{e}_2, \widetilde{e}_3]=0,~~~\alpha+\delta\neq 0,~~~\alpha\gamma-\beta\delta=0.
\end{equation}
where ${\widetilde{e}_1,\widetilde{e}_2,\widetilde{e}_3}$ is a pseudo-orthonormal basis, with $\widetilde{e}_3$ timelike.
\begin{lem} The Bott connection $\nabla^B$ of $G_6$ is given by
\begin{align}
&\nabla^B_{\widetilde{e}_1}\widetilde{e}_1=0,~~~\nabla^B_{\widetilde{e}_1}\widetilde{e}_2=0,~~~\nabla^B_{\widetilde{e}_1}\widetilde{e}_3=\delta\widetilde{e}_3,\nonumber\\
&\nabla^B_{\widetilde{e}_2}\widetilde{e}_1=-\alpha\widetilde{e}_2,~~~\nabla^B_{\widetilde{e}_2}\widetilde{e}_2=\alpha\widetilde{e}_1,~~~\nabla^B_{\widetilde{e}_2}\widetilde{e}_3=0,\nonumber\\
&\nabla^B_{\widetilde{e}_3}\widetilde{e}_1=-\gamma\widetilde{e}_2,~~~\nabla^B_{\widetilde{e}_3}\widetilde{e}_2=0,~~~\nabla^B_{\widetilde{e}_3}\widetilde{e}_3=0.
\end{align}
\end{lem}
\begin{lem} The curvature $R^B$ of the Bott connection $\nabla^B$ of $(G_6,g)$ is given by
\begin{align}
&R^B(\widetilde{e}_1,\widetilde{e}_2)\widetilde{e}_1=(\alpha^2+\beta\gamma)\widetilde{e}_2,~~~R^B(\widetilde{e}_1,\widetilde{e}_2)\widetilde{e}_2=-\alpha^2\widetilde{e}_1,~~~R^B(\widetilde{e}_1,\widetilde{e}_2)\widetilde{e}_3=0,\nonumber\\
&R^B(\widetilde{e}_1,\widetilde{e}_3)\widetilde{e}_1=\gamma(\alpha+\delta)\widetilde{e}_2,~~~R^B(\widetilde{e}_1,\widehat{e}_3)\widetilde{e}_2=-\alpha\gamma\widetilde{e}_1,~~~R^B(\widetilde{e}_1,\widetilde{e}_3)\widetilde{e}_3=0,\nonumber\\
&R^B(\widetilde{e}_2,\widetilde{e}_3)\widetilde{e}_1=-\alpha\gamma\widetilde{e}_1,~~~R^B(\widetilde{e}_2,\widetilde{e}_3)\widetilde{e}_2=\alpha\gamma\widetilde{e}_2,~~~R^B(\widetilde{e}_2,\widetilde{e}_3)\widetilde{e}_3=0.
\end{align}
\end{lem}
By (2.3), we have
\begin{align}
&\rho^B(\widetilde{e}_1,\widetilde{e}_1)=-(\alpha^2+\beta\gamma),~~~\rho^B(\widetilde{e}_1,\widetilde{e}_2)=\rho^B(\widetilde{e}_1,\widetilde{e}_3)=0,\nonumber\\
&\rho^B(\widetilde{e}_2,\widetilde{e}_1)=0,~~~\rho^B(\widetilde{e}_2,\widetilde{e}_2)=-\alpha^2,~~~\rho^B(\widetilde{e}_2,\widetilde{e}_3)=0,\nonumber\\
&\rho^B(\widetilde{e}_3,\widetilde{e}_1)=\rho^B(\widetilde{e}_3,\widetilde{e}_2)=\rho^B(\widetilde{e}_3,\widetilde{e}_3)=0.
\end{align}
Then,
\begin{align}
&\widetilde{\rho}^B(\widetilde{e}_1,\widetilde{e}_1)=-(\alpha^2+\beta\gamma),~~~\widetilde{\rho}^B(\widetilde{e}_1,\widetilde{e}_2)=\widetilde{\rho}^B(\widetilde{e}_1,\widetilde{e}_3)=0,\nonumber\\
&\widetilde{\rho}^B(\widetilde{e}_2,\widetilde{e}_2)=-\alpha^2,~~~\widetilde{\rho}^B(\widetilde{e}_2,\widetilde{e}_3)=0,~~~\widetilde{\rho}^B(\widetilde{e}_3,\widetilde{e}_3)=0.
\end{align}
By (2.5), we have
\begin{align}
&(\nabla^B_{\widetilde{e}_1}\widetilde{\rho}^B)(\widetilde{e}_2,\widetilde{e}_1)=(\nabla^B_{\widetilde{e}_2}\widetilde{\rho}^B)(\widetilde{e}_1,\widetilde{e}_1)=(\nabla^B_{\widetilde{e}_1}\widetilde{\rho}^B)(\widetilde{e}_2,\widetilde{e}_2)=0,\nonumber\\
&(\nabla^B_{\widetilde{e}_2}\widetilde{\rho}^B)(\widetilde{e}_1,\widetilde{e}_2)=\alpha\beta\gamma,~~~(\nabla^B_{\widetilde{e}_1}\widetilde{\rho}^B)(\widetilde{e}_2,\widetilde{e}_3)=0,~~~(\nabla^B_{\widetilde{e}_2}\widetilde{\rho}^B)(\widetilde{e}_1,\widetilde{e}_3)=0,\nonumber\\
&(\nabla^B_{\widetilde{e}_1}\widetilde{\rho}^B)(\widetilde{e}_3,\widetilde{e}_1)=(\nabla^B_{\widetilde{e}_3}\widetilde{\rho}^B)(\widetilde{e}_1,\widetilde{e}_1)=(\nabla^B_{\widetilde{e}_1}\widetilde{\rho}^B)(\widetilde{e}_3,\widetilde{e}_2)=0,\nonumber\\
&(\nabla^B_{\widetilde{e}_3}\widetilde{\rho}^B)(\widetilde{e}_1,\widetilde{e}_2)=-\alpha^2\gamma,~~~(\nabla^B_{\widetilde{e}_1}\widetilde{\rho}^B)(\widetilde{e}_3,\widetilde{e}_3)=0,~~~(\nabla^B_{\widetilde{e}_3}\widetilde{\rho}^B)(\widetilde{e}_1,\widetilde{e}_3)=0,\nonumber\\
&(\nabla^B_{\widetilde{e}_2}\widetilde{\rho}^B)(\widetilde{e}_3,\widetilde{e}_1)=0,~~~(\nabla^B_{\widetilde{e}_3}\widetilde{\rho}^B)(\widetilde{e}_2,\widetilde{e}_1)=-\alpha^2\gamma,~~~(\nabla^B_{\widetilde{e}_2}\widetilde{\rho}^B)(\widetilde{e}_3,\widetilde{e}_2)=0,\nonumber\\
&(\nabla^B_{\widetilde{e}_3}\widetilde{\rho}^B)(\widetilde{e}_2,\widetilde{e}_2)=(\nabla^B_{\widetilde{e}_2}\widetilde{\rho}^B)(\widetilde{e}_3,\widetilde{e}_3)=(\nabla^B_{\widetilde{e}_3}\widetilde{\rho}^B)(\widetilde{e}_2,\widetilde{e}_3)=0.
\end{align}
Then, if $\widetilde{\rho}^B$  is a Codazzi tensor on $(G_6,\nabla^B)$, by (2.6) and (2.7), we have the following two equations:
\begin{eqnarray}
       \begin{cases}
        \alpha\beta\gamma=0 \\[2pt]
       \alpha^2\gamma=0\\[2pt]
       \end{cases}
\end{eqnarray}
By solving (2.47), we get\\
\begin{thm}
 $\widetilde{\rho}^B$ is a Codazzi tensor on $(G_6,\nabla^B)$ if and only if\\
\begin{eqnarray*}
&&(1)\alpha=\beta=0,~~~ \delta\neq 0;\nonumber\\
&&(2)\alpha\neq0,~~~\gamma=\beta\delta=0.\nonumber\\
\end{eqnarray*}
\end{thm}
\vskip 0.5 true cm
\noindent{\bf 2.7 Codazzi tensors  of $G_7$}\\
\vskip 0.5 true cm\indent By \cite{W}, we have the following Lie algebra of $G_7$ satisfies
\begin{equation}
[\widetilde{e}_1, \widetilde{e}_2]=-\alpha\widetilde{e}_1-\beta\widetilde{e}_2-\beta\widetilde{e}_3,~~~[\widetilde{e}_1, \widetilde{e}_3]=\alpha\widetilde{e}_1+\beta\widetilde{e}_2+\beta\widetilde{e}_3,~~~[\widetilde{e}_2, \widetilde{e}_3]=\gamma\widetilde{e}_1+\delta\widetilde{e}_2+\delta\widetilde{e}_3,~~~\alpha+\delta\neq 0,~~~\alpha\gamma=0.
\end{equation}
where ${\widetilde{e}_1,\widetilde{e}_2,\widetilde{e}_3}$ is a pseudo-orthonormal basis, with $\widetilde{e}_3$ timelike.
\begin{lem} The Bott connection $\nabla^B$ of $G_7$ is given by
\begin{align}
&\nabla^B_{\widetilde{e}_1}\widetilde{e}_1=\alpha\widetilde{e}_2,~~~\nabla^B_{\widetilde{e}_1}\widetilde{e}_2=-\alpha\widetilde{e}_1,~~~\nabla^B_{\widetilde{e}_1}\widetilde{e}_3=\beta\widetilde{e}_3,\nonumber\\
&\nabla^B_{\widetilde{e}_2}\widetilde{e}_1=\beta\widetilde{e}_2,~~~\nabla^B_{\widetilde{e}_2}\widetilde{e}_2=-\beta\widetilde{e}_1,~~~\nabla^B_{\widetilde{e}_2}\widetilde{e}_3=\delta\widetilde{e}_3,\nonumber\\
&\nabla^B_{\widetilde{e}_3}\widetilde{e}_1=-\alpha\widetilde{e}_1-\beta\widetilde{e}_2,~~~\nabla^B_{\widetilde{e}_3}\widetilde{e}_2=-\gamma\widetilde{e}_1-\delta\widetilde{e}_2,~~~\nabla^B_{\widetilde{e}_3}\widetilde{e}_3=0.
\end{align}
\end{lem}
\begin{lem} The curvature $R^B$ of the Bott connection $\nabla^B$ of $(G_7,g)$ is given by
\begin{align}
&R^B(\widetilde{e}_1,\widetilde{e}_2)\widetilde{e}_1=-\alpha\beta\widetilde{e}_1+\alpha^2\widetilde{e}_2,~~~R^B(\widetilde{e}_1,\widetilde{e}_2)\widetilde{e}_2=-(\alpha^2+\beta^2+\beta\gamma)\widetilde{e}_1-\beta\delta\widetilde{e}_2,~~~R^B(\widetilde{e}_1,\widetilde{e}_2)\widetilde{e}_3=\beta(\alpha-\delta)\widetilde{e}_3,\nonumber\\
&R^B(\widetilde{e}_1,\widetilde{e}_3)\widetilde{e}_1=\alpha(2\beta+\gamma)\widetilde{e}_1+(\alpha\delta-2\alpha^2)\widetilde{e}_2,~~~R^B(\widetilde{e}_1,\widetilde{e}_3)\widetilde{e}_2=(\alpha\delta+\beta^2+\beta\gamma)\widetilde{e}_1+(\beta\delta-\alpha\beta-\alpha\gamma)\widetilde{e}_2,\nonumber\\
&R^B(\widetilde{e}_1,\widetilde{e}_3)\widetilde{e}_3=-\beta(\alpha+\delta)\widetilde{e}_3,~~~R^B(\widetilde{e}_2,\widetilde{e}_3)\widetilde{e}_1=(\beta^2+\beta\gamma+\alpha\delta)\widetilde{e}_1+(\beta\delta-\alpha\beta-\alpha\gamma)\widetilde{e}_2,\nonumber\\
&R^B(\widetilde{e}_2,\widetilde{e}_3)\widetilde{e}_2=(2\beta\delta+\delta\gamma+\alpha\gamma-\alpha\beta)\widetilde{e}_1+(\delta^2-\beta^2-\beta\gamma)\widetilde{e}_2,~~~R^B(\widetilde{e}_2,\widetilde{e}_3)\widetilde{e}_3=-(\beta\gamma+\delta^2)\widetilde{e}_3.
\end{align}
\end{lem}
By (2.3), we have
\begin{align}
&\rho^B(\widetilde{e}_1,\widetilde{e}_1)=-\alpha^2,~~~\rho^B(\widetilde{e}_1,\widetilde{e}_2)=\beta\delta,~~~\rho^B(\widetilde{e}_1,\widetilde{e}_3)=\beta(\alpha+\delta),\nonumber\\
&\rho^B(\widetilde{e}_2,\widetilde{e}_1)=-\alpha\beta,~~~\rho^B(\widetilde{e}_2,\widetilde{e}_2)=-(\alpha^2+\beta^2+\beta\gamma),~~~\rho^B(\widetilde{e}_2,\widetilde{e}_3)=(\beta\gamma+\delta^2),\nonumber\\
&\rho^B(\widetilde{e}_3,\widetilde{e}_1)=\beta(\alpha+\delta),~~~\rho^B(\widetilde{e}_3,\widetilde{e}_2)=\delta(\alpha+\delta),~~~\rho^B(\widetilde{e}_3,\widetilde{e}_3)=0.
\end{align}
Then,
\begin{align}
&\widetilde{\rho}^B(\widetilde{e}_1,\widetilde{e}_1)=-\alpha^2,~~~\widetilde{\rho}^B(\widetilde{e}_1,\widetilde{e}_2)=\frac{\beta(\delta-\alpha)}{2},~~~\widetilde{\rho}^B(\widetilde{e}_1,\widetilde{e}_3)=\delta(\alpha+\delta),\nonumber\\
&\widetilde{\rho}^B(\widetilde{e}_2,\widetilde{e}_2)=-(\alpha^2+\beta^2+\beta\gamma),~~~\widetilde{\rho}^B(\widetilde{e}_2,\widetilde{e}_3)=\delta^2+\frac{\beta\gamma+\alpha\delta}{2},~~~\widetilde{\rho}^B(\widetilde{e}_3,\widetilde{e}_3)=0.
\end{align}
By (2.5), we have
\begin{align}
&(\nabla^B_{\widetilde{e}_1}\widetilde{\rho}^B)(\widetilde{e}_2,\widetilde{e}_1)=\alpha(\beta^2+\beta\gamma),~~~(\nabla^B_{\widetilde{e}_2}\widetilde{\rho}^B)(\widetilde{e}_1,\widetilde{e}_1)=\beta^2(\alpha-\delta),~~~(\nabla^B_{\widetilde{e}_1}\widetilde{\rho}^B)(\widetilde{e}_2,\widetilde{e}_2)=\alpha\beta(\delta-\alpha),\nonumber\\
&(\nabla^B_{\widetilde{e}_2}\widetilde{\rho}^B)(\widetilde{e}_1,\widetilde{e}_2)=\beta^2(\beta+\gamma),~~~(\nabla^B_{\widetilde{e}_1}\widetilde{\rho}^B)(\widetilde{e}_2,\widetilde{e}_3)=\alpha^2\beta+\frac{\alpha\beta\delta-\beta^2\gamma}{2}-\beta\delta^2,\nonumber\\
&(\nabla^B_{\widetilde{e}_2}\widetilde{\rho}^B)(\widetilde{e}_1,\widetilde{e}_3)=-(2\beta\delta^2+\frac{\beta^2\gamma+3\alpha\beta\delta}{2}),~~~(\nabla^B_{\widetilde{e}_1}\widetilde{\rho}^B)(\widetilde{e}_3,\widetilde{e}_1)=-(\alpha\beta^2+\beta^2\delta+\alpha\delta^2+\frac{\alpha\beta\gamma+\alpha^2\delta}{2}),\nonumber\\
&(\nabla^B_{\widetilde{e}_3}\widetilde{\rho}^B)(\widetilde{e}_1,\widetilde{e}_1)=\beta^2\delta-\alpha\beta^2-2\alpha^3,~~~(\nabla^B_{\widetilde{e}_1}\widetilde{\rho}^B)(\widetilde{e}_3,\widetilde{e}_2)=\alpha^2\beta+\frac{\alpha\beta\delta}{2}-\beta\delta^2-\frac{\beta^2\gamma}{2},\nonumber\\
&(\nabla^B_{\widetilde{e}_3}\widetilde{\rho}^B)(\widetilde{e}_1,\widetilde{e}_2)=\frac{\beta(\delta^2-3\alpha^2)}{2}-\beta^3-\beta^2\gamma-\alpha^2\gamma,~~~(\nabla^B_{\widetilde{e}_1}\widetilde{\rho}^B)(\widetilde{e}_3,\widetilde{e}_3)=0,\nonumber\\
&(\nabla^B_{\widetilde{e}_3}\widetilde{\rho}^B)(\widetilde{e}_1,\widetilde{e}_3)=\alpha^2\beta+\frac{3\alpha\beta\delta}{2}+\beta\delta^2+\frac{\beta^2\gamma}{2},~~~(\nabla^B_{\widetilde{e}_2}\widetilde{\rho}^B)(\widetilde{e}_3,\widetilde{e}_1)=-(2\beta\delta^2+\frac{3\alpha\beta\delta+\beta^2\gamma}{2}),\nonumber\\
&(\nabla^B_{\widetilde{e}_3}\widetilde{\rho}^B)(\widetilde{e}_2,\widetilde{e}_1)=\frac{\beta\delta^2-3\alpha^2\beta}{2}-\alpha^2\gamma-\beta^3-\beta^2\gamma,~~~(\nabla^B_{\widetilde{e}_2}\widetilde{\rho}^B)(\widetilde{e}_3,\widetilde{e}_2)=\beta^2(\alpha+\delta)-\delta^3-\frac{\beta\delta\gamma+\alpha\delta^2}{2},\nonumber\\
&(\nabla^B_{\widetilde{e}_3}\widetilde{\rho}^B)(\widetilde{e}_2,\widetilde{e}_2)=-(\alpha\beta\gamma+\beta\delta\gamma+2\alpha^2\delta+2\beta^2\delta),~~~(\nabla^B_{\widetilde{e}_2}\widetilde{\rho}^B)(\widetilde{e}_3,\widetilde{e}_3)=0,~~~(\nabla^B_{\widetilde{e}_3}\widetilde{\rho}^B)(\widetilde{e}_2,\widetilde{e}_3)=\alpha\beta\gamma+\delta^3+\frac{3\beta\delta\gamma+\alpha\delta^2}{2}.
\end{align}
Then, if $\widetilde{\rho}^B$  is a Codazzi tensor on $(G_7,\nabla^B)$, by (2.6) and (2.7), we have the following nine equations:
\begin{eqnarray}
       \begin{cases}
        \beta(\alpha\gamma+\beta\delta)=0 \\[2pt]
       \beta(\alpha\delta-\alpha^2-\beta^2-\beta\gamma)=0\\[2pt]
       \beta(\alpha+\delta)^2=0 \\[2pt]
       2\alpha^3-2\beta^2\delta-\alpha\delta^2-\frac{\alpha\beta\gamma+\alpha^2\delta}{2}=0 \\[2pt]
       \frac{5\alpha^2\beta+\alpha\beta\delta+\beta^2\gamma-3\beta\delta^2}{2}+\alpha^2\gamma+\beta^3=0 \\[2pt]
       \beta(\alpha^2+3\alpha\delta+\delta^2+\frac{\beta\gamma}{2})=0 \\[2pt]
       \frac{3\alpha^2\beta-3\alpha\beta\delta+\beta^2\gamma-5\beta\delta^2}{2}+\beta^3+\alpha^2\gamma=0   \\[2pt]    \frac{\beta\delta\gamma-\alpha\delta^2}{2}+\alpha\beta\gamma+2\alpha^2\delta+\alpha\beta^2+3\beta^2\delta-\delta^3=0\\[2pt]
       \alpha\beta\gamma+\delta^3+\frac{3\beta\delta\gamma+\alpha\delta^2}{2}=0 \\[2pt]
       \end{cases}
\end{eqnarray}
By solving (2.54), we get $\alpha=\delta=0$, there is a contradiction. So\\
\begin{thm}
 $\widetilde{\rho}^B$  is not a Codazzi tensor on $(G_7,\nabla^B)$.
\end{thm}
\section{Quasi-statistical structure associated to Bott connections on three-dimensional Lorentzian Lie groups}
The torsion tensor of $(G_i,g,\nabla^B)$ is defined by
\begin{equation}
T^B(X,Y)=\nabla^B_XY-\nabla^B_YX-[X,Y].
\end{equation}
Then we have
\begin{defn} \cite{A} Let $M$ be a smooth manifold endowed with a linear connection $\nabla$, and a tensor fields $\omega$. Then $(M,\nabla,\omega)$ is called a quasi-statistical structure, if it satisfies
\begin{equation}
\widetilde{f}(X,Y,Z)=(\nabla_X\omega)(Y,Z)-(\nabla_Y\omega)(X,Z)+\omega(T(X,Y),Z)=0,
\end{equation}
where $\widetilde{f}$ is $C^\infty(M)$-linear for $X,Y,Z$, and $\widetilde{f}(X,Y,Z)=-\widetilde{f}(Y,X,Z)$.
\end{defn}
Then we have $(G_i,\nabla^B,\omega)$ is a quasi-statistical structure if and only if the following nine equations hold:
\begin{eqnarray}
       \begin{cases}
\widetilde{f}(\widetilde{e}_1,\widetilde{e}_2,\widetilde{e}_j)=0\\[2pt]
\widetilde{f}(\widetilde{e}_1,\widetilde{e}_3,\widetilde{e}_j)=0\\[2pt]
\widetilde{f}(\widetilde{e}_2,\widetilde{e}_3,\widetilde{e}_j)=0\\[2pt]
 \end{cases}
\end{eqnarray}
where $1\leq j\leq3$.\\
For $(G_1,\nabla^B)$, we have
\begin{align}
T^B(\widetilde{e}_1,\widetilde{e}_2)=\beta \widetilde{e}_3,~~~T^B(\widetilde{e}_1,\widetilde{e}_3)=T^B(\widetilde{e}_2,\widetilde{e}_3)=0.
\end{align}
\begin{align}
&\widetilde{\rho}^B(T^B(\widetilde{e}_1,\widetilde{e}_2),\widetilde{e}_1)=-\frac{\alpha\beta^2}{2},~~~\widetilde{\rho}^B(T^B(\widetilde{e}_1,\widetilde{e}_2),\widetilde{e}_2)=\frac{\alpha\beta^2}{2},~~~\widetilde{\rho}^B(T^B(\widetilde{e}_1,\widetilde{e}_2),\widetilde{e}_3)=0,\\\nonumber
&\widetilde{\rho}^B(T^B(\widetilde{e}_1,\widetilde{e}_3),\widetilde{e}_j)=\widetilde{\rho}^B(T^B(\widetilde{e}_2,\widetilde{e}_3),\widetilde{e}_j)=0,
\end{align}
where $1\leq j\leq3$.\\
Then, if $(G_1,\nabla^B,\widetilde{\rho}^B)$  is a quasi-statistical structure, by (3.2) and (3.3), we have the following three equations:
\begin{eqnarray}
        \begin{cases}
        2\alpha^2\beta=0 \\[2pt]
      \frac{3\alpha^3}{2}=0\\[2pt]
       \frac{\alpha}{2}(\alpha^2-\beta^2)=0\\[2pt]
       \end{cases}
\end{eqnarray}
By solving (3.6), we get $\alpha=0$, there is a contradiction. So\\
\begin{thm}
 $(G_1,\nabla^B,\widetilde{\rho}^B)$ is not a quasi-statistical structure.
\end{thm}
For $(G_2,\nabla^B)$, we have
\begin{align}
T^B(\widetilde{e}_1,\widetilde{e}_2)=\beta \widetilde{e}_3,~~~T^B(\widetilde{e}_1,\widetilde{e}_3)=T^B(\widetilde{e}_2,\widetilde{e}_3)=0.
\end{align}
\begin{align}
&\widetilde{\rho}^B(T^B(\widetilde{e}_1,\widetilde{e}_2),\widetilde{e}_1)=0,~~~\widetilde{\rho}^B(T^B(\widetilde{e}_1,\widetilde{e}_2),\widetilde{e}_2)=-\frac{\alpha\beta\gamma}{2},~~~\widetilde{\rho}^B(T^B(\widetilde{e}_1,\widetilde{e}_2),\widetilde{e}_3)=0,\\\nonumber
&\widetilde{\rho}^B(T^B(\widetilde{e}_1,\widetilde{e}_3),\widetilde{e}_j)=\widetilde{\rho}^B(T^B(\widetilde{e}_2,\widetilde{e}_3),\widetilde{e}_j)=0,
\end{align}
where $1\leq j\leq3$.\\
Then, if $(G_2,\nabla^B,\widetilde{\rho}^B)$ is a quasi-statistical structure, by (3.2) and (3.3), we have the following three equations:
\begin{eqnarray}
        \begin{cases}
        \gamma(\frac{\alpha\beta}{2}-\beta^2)=0 \\[2pt]
      \frac{\alpha\beta\gamma}{2}=0\\[2pt]
       \gamma^2(\frac{\alpha}{2}-\beta)=0\\[2pt]
       \end{cases}
\end{eqnarray}
By solving (3.9), we get\\
\begin{thm}
 $(G_2,\nabla^B,\widetilde{\rho}^B)$ is a quasi-statistical structure if and only if
$\alpha=\beta=0,~~~\gamma\neq0$.
\end{thm}
For $(G_3,\nabla^B)$, we have
\begin{align}
T^B(\widetilde{e}_1,\widetilde{e}_2)=\gamma \widetilde{e}_3,~~~T^B(\widetilde{e}_1,\widetilde{e}_3)=T^B(\widetilde{e}_2,\widetilde{e}_3)=0.
\end{align}
\begin{align}
&\widetilde{\rho}^B(T^B(\widetilde{e}_1,\widetilde{e}_2),\widetilde{e}_j)=\widetilde{\rho}^B(T^B(\widetilde{e}_1,\widetilde{e}_3),\widetilde{e}_j)=\widetilde{\rho}^B(T^B(\widetilde{e}_2,\widetilde{e}_3),\widetilde{e}_j)=0,
\end{align}
where $1\leq j\leq3$.\\
Similarly, we can get
\begin{thm}
$(G_3,\nabla^B,\widetilde{\rho}^B)$ is a quasi-statistical structure
\end{thm}
For $(G_4,\nabla^B)$, we have
\begin{align}
T^B(\widetilde{e}_1,\widetilde{e}_2)=(\beta-2\eta) \widetilde{e}_3,~~~T^B(\widetilde{e}_1,\widetilde{e}_3)=T^B(\widetilde{e}_2,\widetilde{e}_3)=0.
\end{align}
\begin{align}
&\widetilde{\rho}^B(T^B(\widetilde{e}_1,\widetilde{e}_2),\widetilde{e}_1)=0,~~~\widetilde{\rho}^B(T^B(\widetilde{e}_1,\widetilde{e}_2),\widetilde{e}_2)=\frac{\alpha(\beta-2\eta)}{2},~~~\widetilde{\rho}^B(T^B(\widetilde{e}_1,\widetilde{e}_2),\widetilde{e}_3)=0,\\\nonumber
&\widetilde{\rho}^B(T^B(\widetilde{e}_1,\widetilde{e}_3),\widetilde{e}_j)=\widetilde{\rho}^B(T^B(\widetilde{e}_2,\widetilde{e}_3),\widetilde{e}_j)=0,
\end{align}
where $1\leq j\leq3$.\\
Then, if $(G_4,\nabla^B,\widetilde{\rho}^B)$ is a quasi-statistical structure, by (3.2) and (3.3), we have the following three equations:
\begin{eqnarray}
        \begin{cases}
        \beta^2-2\beta\eta+\alpha\eta-\frac{\alpha\beta}{2}=0 \\[2pt]
      \frac{\alpha}{2}-\beta=0\\[2pt]
       \frac{\alpha\beta}{2}=0\\[2pt]
       \end{cases}
\end{eqnarray}
By solving (3.14), we get\\
\begin{thm}
 $(G_4,\nabla^B,\widetilde{\rho}^B)$ is a quasi-statistical structure if and only if
$\alpha=\beta=0$.
\end{thm}
For $(G_5,\nabla^B)$, we have
\begin{align}
T^B(\widetilde{e}_1,\widetilde{e}_2)=T^B(\widetilde{e}_1,\widetilde{e}_3)=T^B(\widetilde{e}_2,\widetilde{e}_3)=0.
\end{align}
\begin{align}
&\widetilde{\rho}^B(T^B(\widetilde{e}_1,\widetilde{e}_2),\widetilde{e}_j)=\widetilde{\rho}^B(T^B(\widetilde{e}_1,\widetilde{e}_3),\widetilde{e}_j)=\widetilde{\rho}^B(T^B(\widetilde{e}_2,\widetilde{e}_3),\widetilde{e}_j)=0,
\end{align}
where $1\leq j\leq3$.\\
Similarly, we can get
\begin{thm}
$(G_5,\nabla^B,\widetilde{\rho}^B)$ is a quasi-statistical structure.
\end{thm}
For $(G_6,\nabla^B)$, we have
\begin{align}
T^B(\widetilde{e}_1,\widetilde{e}_2)=-\beta \widetilde{e}_3,~~~T^B(\widetilde{e}_1,\widetilde{e}_3)=T^B(\widetilde{e}_2,\widetilde{e}_3)=0.
\end{align}
\begin{align}
&\widetilde{\rho}^B(T^B(\widetilde{e}_1,\widetilde{e}_2),\widetilde{e}_1)=\widetilde{\rho}^B(T^B(\widetilde{e}_1,\widetilde{e}_3),\widetilde{e}_j)=\widetilde{\rho}^B(T^B(\widetilde{e}_2,\widetilde{e}_3),\widetilde{e}_j)=0,
\end{align}
where $1\leq j\leq3$.\\
Then, if $(G_6,\nabla^B,\widetilde{\rho}^B)$ is a quasi-statistical structure, by (3.2) and (3.3), we have the following two equations:
\begin{eqnarray}
       \begin{cases}
        \alpha\beta\gamma=0 \\[2pt]
       \alpha^2\gamma=0\\[2pt]
       \end{cases}
\end{eqnarray}
By solving (3.19), we get\\
\begin{thm}
 $(G_6,\nabla^B,\widetilde{\rho}^B)$ is a quasi-statistical structure if and only if\\
\begin{eqnarray*}
&&(1)\alpha=\beta=0,~~~ \delta\neq 0;\nonumber\\
&&(2)\alpha\neq0,~~~\gamma=\beta\delta=0.\nonumber\\
\end{eqnarray*}
\end{thm}
For $(G_7,\nabla^B)$, we have
\begin{align}
T^B(\widetilde{e}_1,\widetilde{e}_2)=\beta \widetilde{e}_3,~~~T^B(\widetilde{e}_1,\widetilde{e}_3)=T^B(\widetilde{e}_2,\widetilde{e}_3)=0.
\end{align}
\begin{align}
&\widetilde{\rho}^B(T^B(\widetilde{e}_1,\widetilde{e}_2),\widetilde{e}_1)=\beta^2(\alpha+\delta),~~~\widetilde{\rho}^B(T^B(\widetilde{e}_1,\widetilde{e}_2),\widetilde{e}_2)=\beta\delta^2+\frac{\alpha\beta\delta+\beta^2\gamma}{2},~~~\widetilde{\rho}^B(T^B(\widetilde{e}_1,\widetilde{e}_2),\widetilde{e}_3)=0,\\\nonumber
&\widetilde{\rho}^B(T^B(\widetilde{e}_1,\widetilde{e}_3),\widetilde{e}_j)=\widetilde{\rho}^B(T^B(\widetilde{e}_2,\widetilde{e}_3),\widetilde{e}_j)=0,
\end{align}
where $1\leq j\leq3$.\\
Then, if $(G_7,\nabla^B,\widetilde{\rho}^B)$ is a quasi-statistical structure, by (3.2) and (3.3), we have the following nine equations:
\begin{eqnarray}
       \begin{cases}
        \beta(\alpha\gamma+\alpha\beta+2\beta\delta)=0 \\[2pt]
       \beta(\alpha\delta-\alpha^2-\beta^2-\beta\gamma+\delta^2+\frac{\beta\gamma+\alpha\delta}{2})=0\\[2pt]
       \beta(\alpha+\delta)^2=0 \\[2pt]
       2\alpha^3-2\beta^2\delta-\alpha\delta^2-\frac{\alpha\beta\gamma+\alpha^2\delta}{2}=0 \\[2pt]
       \frac{5\alpha^2\beta+\alpha\beta\delta+\beta^2\gamma-3\beta\delta^2}{2}+\alpha^2\gamma+\beta^3=0 \\[2pt]
       \beta(\alpha^2+3\alpha\delta+\delta^2+\frac{\beta\gamma}{2})=0 \\[2pt]
       \frac{3\alpha^2\beta-3\alpha\beta\delta+\beta^2\gamma-5\beta\delta^2}{2}+\beta^3+\alpha^2\gamma=0   \\[2pt]    \frac{\beta\delta\gamma-\alpha\delta^2}{2}+\alpha\beta\gamma+2\alpha^2\delta+\alpha\beta^2+3\beta^2\delta-\delta^3=0\\[2pt]
       \alpha\beta\gamma+\delta^3+\frac{3\beta\delta\gamma+\alpha\delta^2}{2}=0 \\[2pt]
       \end{cases}
\end{eqnarray}
By solving (3.22), we get $\alpha=\delta=0$, there is a contradiction. So\\
\begin{thm}
 $(G_7,\nabla^B,\widetilde{\rho}^B)$ is not a quasi-statistical structure.
\end{thm}
\section{Codazzi tensors associated to canonical connections and Kobayashi-Nomizu connections on three-dimensional Lorentzian Lie groups}
By \cite{FE}, we define canonical connections and Kobayashi-Nomizu connections as follows:
\begin{eqnarray}
\nabla^c_XY=\nabla^L_XY-\frac{1}{2}(\nabla_XJ)JY,
\end{eqnarray}
\begin{eqnarray}
\nabla^k_XY=\nabla^c_XY-\frac{1}{4}[(\nabla_YJ)JX-(\nabla_{JY}J)X],
\end{eqnarray}
where $J$ is a product structure on $\{G_i\}_{i=1,2\cdot\cdot\cdot,7}$ by $J\widetilde{e}_1=\widetilde{e}_1,J\widetilde{e}_2=\widetilde{e}_2,J\widetilde{e}_3=-\widetilde{e}_3$.\\
\vskip 0.5 true cm
\noindent{\bf 4.1 Codazzi tensors of $G_1$}\\
\vskip 0.5 true cm
\begin{lem} (\cite{wy}) The canonical connection $\nabla^c$ of $G_1$ is given by
\begin{align}
&\nabla^c_{\widetilde{e}_1}\widetilde{e}_1=-\alpha\widetilde{e}_2,~~~\nabla^c_{\widetilde{e}_1}\widetilde{e}_2=\alpha\widetilde{e}_1,~~~\nabla^c_{\widetilde{e}_1}\widetilde{e}_3=0,\nonumber\\
&\nabla^c_{\widetilde{e}_2}\widetilde{e}_1=\nabla^c_{\widetilde{e}_2}\widetilde{e}_2=\nabla^c_{\widetilde{e}_2}\widetilde{e}_3=0,\nonumber\\
&\nabla^c_{\widetilde{e}_3}\widetilde{e}_1=\frac{\beta}{2}\widetilde{e}_2,~~~\nabla^c_{\widetilde{e}_3}\widetilde{e}_2=-\frac{\beta}{2}\widetilde{e}_1,~~~\nabla^c_{\widetilde{e}_3}\widetilde{e}_3=0.
\end{align}
\end{lem}
Then,
\begin{align}
&\widetilde{\rho}^c(\widetilde{e}_1,\widetilde{e}_1)=-(\alpha^2+\frac{\beta^2}{2}),~~~\widetilde{\rho}^c(\widetilde{e}_1,\widetilde{e}_2)=0,~~~\widetilde{\rho}^c(\widetilde{e}_1,\widetilde{e}_3)=\frac{\alpha\beta}{4},\nonumber\\
&\widetilde{\rho}^c(\widetilde{e}_2,\widetilde{e}_2)=-(\alpha^2+\frac{\beta^2}{2}),~~~\widetilde{\rho}^c(\widetilde{e}_2,\widetilde{e}_3)=\frac{\alpha^2}{2},~~~\widetilde{\rho}^c(\widetilde{e}_3,\widetilde{e}_3)=0.
\end{align}
By (2.5), we have
\begin{align}
&(\nabla^c_{\widetilde{e}_1}\widetilde{\rho}^c)(\widetilde{e}_2,\widetilde{e}_1)=(\nabla^c_{\widetilde{e}_2}\widetilde{\rho}^c)(\widetilde{e}_1,\widetilde{e}_1)=(\nabla^c_{\widetilde{e}_1}\widetilde{\rho}^c)(\widetilde{e}_2,\widetilde{e}_2)=0,\nonumber\\
&(\nabla^c_{\widetilde{e}_2}\widetilde{\rho}^c)(\widetilde{e}_1,\widetilde{e}_2)=0,~~~(\nabla^c_{\widetilde{e}_1}\widetilde{\rho}^c)(\widetilde{e}_2,\widetilde{e}_3)=-\frac{\alpha\beta}{4},~~~(\nabla^c_{\widetilde{e}_2}\widetilde{\rho}^c)(\widetilde{e}_1,\widetilde{e}_3)=-\frac{\alpha\beta}{4},\nonumber\\
&(\nabla^c_{\widetilde{e}_1}\widetilde{\rho}^c)(\widetilde{e}_3,\widetilde{e}_1)=\frac{\alpha^3}{2},~~~(\nabla^c_{\widetilde{e}_3}\widetilde{\rho}^c)(\widetilde{e}_1,\widetilde{e}_1)=0,~~~(\nabla^c_{\widetilde{e}_1}\widetilde{\rho}^c)(\widetilde{e}_3,\widetilde{e}_2)=\frac{\alpha^2\beta}{2},\nonumber\\
&(\nabla^c_{\widetilde{e}_3}\widetilde{\rho}^c)(\widetilde{e}_1,\widetilde{e}_2)=(\nabla^c_{\widetilde{e}_1}\widetilde{\rho}^c)(\widetilde{e}_3,\widetilde{e}_3)=0,~~~(\nabla^c_{\widetilde{e}_3}\widetilde{\rho}^c)(\widetilde{e}_1,\widetilde{e}_3)=-\frac{\alpha^2\beta}{4},\nonumber\\
&(\nabla^c_{\widetilde{e}_2}\widetilde{\rho}^c)(\widetilde{e}_3,\widetilde{e}_1)=(\nabla^c_{\widetilde{e}_3}\widetilde{\rho}^c)(\widetilde{e}_2,\widetilde{e}_1)=0,~~~(\nabla^c_{\widetilde{e}_2}\widetilde{\rho}^c)(\widetilde{e}_3,\widetilde{e}_2)=0,\nonumber\\
&(\nabla^c_{\widetilde{e}_3}\widetilde{\rho}^c)(\widetilde{e}_2,\widetilde{e}_2)=-2\alpha^2,~~~(\nabla^c_{\widetilde{e}_2}\widetilde{\rho}^c)(\widetilde{e}_3,\widetilde{e}_3)=0,~~~(\nabla^c_{\widetilde{e}_3}\widetilde{\rho}^c)(\widetilde{e}_2,\widetilde{e}_3)=\frac{\alpha\beta^2}{8}.
\end{align}
Then, if $\widetilde{\rho}^c$  is a Codazzi tensor on $(G_1,\nabla^c)$, by (2.6) and (2.7), we have the following two equations:
\begin{eqnarray}
       \begin{cases}
      \frac{3\alpha^3}{2}=0\\[2pt]
       \frac{\alpha^2\beta}{4}=0\\[2pt]
       \end{cases}
\end{eqnarray}
By solving (4.6) , we get $\alpha=0$, there is a contradiction. So\\
\begin{thm}
$\widetilde{\rho}^c$  is not a Codazzi tensor on $(G_1,\nabla^c)$.
\end{thm}
\begin{lem} (\cite{wy}) The Kobayashi-Nomizu connection $\nabla^k$ of $G_1$ is given by
\begin{align}
&\nabla^k_{\widetilde{e}_1}\widetilde{e}_1=-\alpha\widetilde{e}_2,~~~\nabla^k_{\widetilde{e}_1}\widetilde{e}_2=\alpha\widetilde{e}_1,~~~\nabla^k_{\widetilde{e}_1}\widetilde{e}_3=0,\nonumber\\
&\nabla^k_{\widetilde{e}_2}\widetilde{e}_1=0,~~~\nabla^k_{\widetilde{e}_2}\widetilde{e}_2=0,~~~\nabla^k_{\widetilde{e}_2}\widetilde{e}_3=\alpha\widetilde{e}_3,\nonumber\\
&\nabla^k_{\widetilde{e}_3}\widetilde{e}_1=\alpha\widetilde{e}_1+\beta\widetilde{e}_2,~~~\nabla^k_{\widetilde{e}_3}\widetilde{e}_2=-\beta\widetilde{e}_1-\alpha\widetilde{e}_2,~~~\nabla^k_{\widetilde{e}_3}\widetilde{e}_3=0.
\end{align}
\end{lem}
Then,
\begin{align}
&\widetilde{\rho}^k(\widetilde{e}_1,\widetilde{e}_1)=-(\alpha^2+\beta^2),~~~\widetilde{\rho}^k(\widetilde{e}_1,\widetilde{e}_2)=\alpha\beta,~~~\widetilde{\rho}^k(\widetilde{e}_1,\widetilde{e}_3)=-\frac{\alpha\beta}{2},\nonumber\\
&\widetilde{\rho}^k(\widetilde{e}_2,\widetilde{e}_2)=-(\alpha^2+\beta^2),~~~\widetilde{\rho}^k(\widetilde{e}_2,\widetilde{e}_3)=\frac{\alpha^2}{2},~~~\widetilde{\rho}^k(\widetilde{e}_3,\widetilde{e}_3)=0.
\end{align}
By (2.5), we have
\begin{align}
&(\nabla^k_{\widetilde{e}_1}\widetilde{\rho}^k)(\widetilde{e}_2,\widetilde{e}_1)=0,~~~(\nabla^k_{\widetilde{e}_2}\widetilde{\rho}^k)(\widetilde{e}_1,\widetilde{e}_1)=0,~~~(\nabla^k_{\widetilde{e}_1}\widetilde{\rho}^k)(\widetilde{e}_2,\widetilde{e}_2)=-2\alpha^2\beta,\nonumber\\
&(\nabla^k_{\widetilde{e}_2}\widetilde{\rho}^k)(\widetilde{e}_1,\widetilde{e}_2)=0,~~~(\nabla^k_{\widetilde{e}_1}\widetilde{\rho}^k)(\widetilde{e}_2,\widetilde{e}_3)=\frac{\alpha^2\beta}{2},~~~(\nabla^k_{\widetilde{e}_2}\widetilde{\rho}^k)(\widetilde{e}_1,\widetilde{e}_3)=\frac{\alpha^2\beta}{2},\nonumber\\
&(\nabla^k_{\widetilde{e}_1}\widetilde{\rho}^k)(\widetilde{e}_3,\widetilde{e}_1)=\frac{\alpha^3}{2},~~~(\nabla^k_{\widetilde{e}_3}\widetilde{\rho}^k)(\widetilde{e}_1,\widetilde{e}_1)=2\alpha^3,~~~(\nabla^k_{\widetilde{e}_1}\widetilde{\rho}^k)(\widetilde{e}_3,\widetilde{e}_2)=\frac{\alpha^2\beta}{2},\nonumber\\
&(\nabla^k_{\widetilde{e}_3}\widetilde{\rho}^k)(\widetilde{e}_1,\widetilde{e}_2)=0,~~~(\nabla^k_{\widetilde{e}_1}\widetilde{\rho}^k)(\widetilde{e}_3,\widetilde{e}_3)=0,~~~(\nabla^k_{\widetilde{e}_3}\widetilde{\rho}^k)(\widetilde{e}_1,\widetilde{e}_3)=0,\nonumber\\
&(\nabla^k_{\widetilde{e}_2}\widetilde{\rho}^k)(\widetilde{e}_3,\widetilde{e}_1)=\frac{\alpha^2\beta}{2},~~~(\nabla^k_{\widetilde{e}_3}\widetilde{\rho}^k)(\widetilde{e}_2,\widetilde{e}_1)=0,~~~(\nabla^k_{\widetilde{e}_2}\widetilde{\rho}^k)(\widetilde{e}_3,\widetilde{e}_2)=-\frac{\alpha^3}{2},\nonumber\\
&(\nabla^k_{\widetilde{e}_3}\widetilde{\rho}^k)(\widetilde{e}_2,\widetilde{e}_2)=-2\alpha^3,~~~(\nabla^k_{\widetilde{e}_2}\widetilde{\rho}^k)(\widetilde{e}_3,\widetilde{e}_3)=0,~~~(\nabla^k_{\widetilde{e}_3}\widetilde{\rho}^k)(\widetilde{e}_2,\widetilde{e}_3)=\frac{\alpha}{2}(\alpha^2-\beta^2).
\end{align}
Then, if $\widetilde{\rho}^k$  is a Codazzi tensor on $(G_1,\nabla^k)$, by (2.6) and (2.7), we have the following three equations:
\begin{eqnarray}
       \begin{cases}
      -2\alpha^2\beta=0\\[2pt]
            \frac{3\alpha^3}{2}=0\\[2pt]
       \frac{\alpha}{2}(\beta^2-\alpha^2)=0\\[2pt]
       \end{cases}
\end{eqnarray}
By solving (4.10), we get $\alpha=0$, there is a contradiction. So\\
\begin{thm}
$\widetilde{\rho}^k$  is not a Codazzi tensor on $(G_1,\nabla^k)$.
\end{thm}
\vskip 0.5 true cm
\noindent{\bf 4.2 Codazzi tensors of $G_2$}\\
\vskip 0.5 true cm
\begin{lem} (\cite{wy}) The canonical connection $\nabla^c$ of $G_2$ is given by
\begin{align}
&\nabla^c_{\widetilde{e}_1}\widetilde{e}_1=0,~~~\nabla^c_{\widetilde{e}_1}\widetilde{e}_2=\alpha\widetilde{e}_1,~~~\nabla^c_{\widetilde{e}_1}\widetilde{e}_3=0,\nonumber\\
&\nabla^c_{\widetilde{e}_2}\widetilde{e}_1=-\gamma\widetilde{e}_2,~~~\nabla^c_{\widetilde{e}_2}\widetilde{e}_2=\gamma\widetilde{e}_1,~~~\nabla^c_{\widetilde{e}_2}\widetilde{e}_3=0,\nonumber\\
&\nabla^c_{\widetilde{e}_3}\widetilde{e}_1=\frac{\alpha}{2}\widetilde{e}_2,~~~\nabla^c_{\widetilde{e}_3}\widetilde{e}_2=-\frac{\alpha}{2}\widetilde{e}_1,~~~\nabla^c_{\widetilde{e}_3}\widetilde{e}_3=0.
\end{align}
\end{lem}
Then,
\begin{align}
&\widetilde{\rho}^c(\widetilde{e}_1,\widetilde{e}_1)=-(\gamma^2+\frac{\alpha\beta}{2}),~~~\widetilde{\rho}^c(\widetilde{e}_1,\widetilde{e}_2)=0,~~~\widetilde{\rho}^c(\widetilde{e}_1,\widetilde{e}_3)=0,\nonumber\\
&\widetilde{\rho}^c(\widetilde{e}_2,\widetilde{e}_2)=-(\gamma^2+\frac{\alpha\beta}{2}),~~~\widetilde{\rho}^c(\widetilde{e}_2,\widetilde{e}_3)=\gamma(\frac{\beta}{2}-\frac{\alpha}{4}),~~~\widetilde{\rho}^c(\widetilde{e}_3,\widetilde{e}_3)=0.
\end{align}
By (2.5), we have
\begin{align}
&(\nabla^c_{\widetilde{e}_1}\widetilde{\rho}^c)(\widetilde{e}_2,\widetilde{e}_1)=(\nabla^c_{\widetilde{e}_2}\widetilde{\rho}^c)(\widetilde{e}_1,\widetilde{e}_1)=(\nabla^c_{\widetilde{e}_1}\widetilde{\rho}^c)(\widetilde{e}_2,\widetilde{e}_2)=0,\nonumber\\
&(\nabla^c_{\widetilde{e}_2}\widetilde{\rho}^c)(\widetilde{e}_1,\widetilde{e}_2)=(\nabla^c_{\widetilde{e}_1}\widetilde{\rho}^c)(\widetilde{e}_2,\widetilde{e}_3)=0,~~~(\nabla^c_{\widetilde{e}_2}\widetilde{\rho}^c)(\widetilde{e}_1,\widetilde{e}_3)=\gamma^2(\frac{\alpha}{4}-\frac{\beta}{2}),\nonumber\\
&(\nabla^c_{\widetilde{e}_1}\widetilde{\rho}^c)(\widetilde{e}_3,\widetilde{e}_1)=(\nabla^c_{\widetilde{e}_3}\widetilde{\rho}^c)(\widetilde{e}_1,\widetilde{e}_1)=(\nabla^c_{\widetilde{e}_1}\widetilde{\rho}^c)(\widetilde{e}_3,\widetilde{e}_2)=0,\nonumber\\
&(\nabla^c_{\widetilde{e}_3}\widetilde{\rho}^c)(\widetilde{e}_1,\widetilde{e}_2)=(\nabla^c_{\widetilde{e}_1}\widetilde{\rho}^c)(\widetilde{e}_3,\widetilde{e}_3)=(\nabla^c_{\widetilde{e}_3}\widetilde{\rho}^c)(\widetilde{e}_1,\widetilde{e}_3)=\frac{\alpha\gamma}{4}(\frac{\alpha}{2}-\beta),\nonumber\\
&(\nabla^c_{\widetilde{e}_2}\widetilde{\rho}^c)(\widetilde{e}_3,\widetilde{e}_1)=\gamma^2(\frac{\beta}{2}-\frac{\alpha}{4}),~~~(\nabla^c_{\widetilde{e}_3}\widetilde{\rho}^c)(\widetilde{e}_2,\widetilde{e}_1)=(\nabla^c_{\widetilde{e}_2}\widetilde{\rho}^c)(\widetilde{e}_3,\widetilde{e}_2)=0,\nonumber\\
&(\nabla^c_{\widetilde{e}_3}\widetilde{\rho}^c)(\widetilde{e}_2,\widetilde{e}_2)=(\nabla^c_{\widetilde{e}_2}\widetilde{\rho}^c)(\widetilde{e}_3,\widetilde{e}_3)=(\nabla^c_{\widetilde{e}_3}\widetilde{\rho}^c)(\widetilde{e}_2,\widetilde{e}_3)=0.
\end{align}
Then, if $\widetilde{\rho}^c$  is a Codazzi tensor on $(G_2,\nabla^c)$, by (2.6) and (2.7), we have the following two equations:
\begin{eqnarray}
       \begin{cases}
      \gamma^2(\frac{\alpha}{4}-\frac{\beta}{2})=0\\[2pt]
      \alpha\gamma(\frac{\beta}{4}-\frac{\alpha}{8})=0\\[2pt]
       \end{cases}
\end{eqnarray}
By solving (4.14), we get\\
\begin{thm}
$\widetilde{\rho}^c$  is a Codazzi tensor on $(G_2,\nabla^c)$ if and only if $\gamma\neq0,~~~\alpha=2\beta$.
\end{thm}
\begin{lem} (\cite{wy}) The Kobayashi-Nomizu connection $\nabla^k$ of $G_2$ is given by
\begin{align}
&\nabla^k_{\widetilde{e}_1}\widetilde{e}_1=0,~~~\nabla^k_{\widetilde{e}_1}\widetilde{e}_2=0,~~~\nabla^k_{\widetilde{e}_1}\widetilde{e}_3=-\gamma\widetilde{e}_3,\nonumber\\
&\nabla^k_{\widetilde{e}_2}\widetilde{e}_1=-\gamma\widetilde{e}_2,~~~\nabla^k_{\widetilde{e}_2}\widetilde{e}_2=\gamma\widetilde{e}_1,~~~\nabla^k_{\widetilde{e}_2}\widetilde{e}_3=0,\nonumber\\
&\nabla^k_{\widetilde{e}_3}\widetilde{e}_1=\beta\widetilde{e}_2,~~~\nabla^k_{\widetilde{e}_3}\widetilde{e}_2=-\alpha\widetilde{e}_1,~~~\nabla^k_{\widetilde{e}_3}\widetilde{e}_3=0.
\end{align}
\end{lem}
Then,
\begin{align}
&\widetilde{\rho}^k(\widetilde{e}_1,\widetilde{e}_1)=-(\gamma^2+\beta^2),~~~\widetilde{\rho}^k(\widetilde{e}_1,\widetilde{e}_2)=0,~~~\widetilde{\rho}^k(\widetilde{e}_1,\widetilde{e}_3)=0,\nonumber\\
&\widetilde{\rho}^k(\widetilde{e}_2,\widetilde{e}_2)=-(\gamma^2+\alpha\beta),~~~\widetilde{\rho}^k(\widetilde{e}_2,\widetilde{e}_3)=-\frac{\alpha\gamma}{2},~~~\widetilde{\rho}^k(\widetilde{e}_3,\widetilde{e}_3)=0.
\end{align}
By (2.5), we have
\begin{align}
&(\nabla^k_{\widetilde{e}_1}\widetilde{\rho}^k)(\widetilde{e}_2,\widetilde{e}_1)=0,~~~(\nabla^k_{\widetilde{e}_2}\widetilde{\rho}^k)(\widetilde{e}_1,\widetilde{e}_1)=0,~~~(\nabla^k_{\widetilde{e}_1}\widetilde{\rho}^k)(\widetilde{e}_2,\widetilde{e}_2)=0,\nonumber\\
&(\nabla^k_{\widetilde{e}_2}\widetilde{\rho}^k)(\widetilde{e}_1,\widetilde{e}_2)=\beta\gamma(\beta-\alpha),~~~(\nabla^k_{\widetilde{e}_1}\widetilde{\rho}^k)(\widetilde{e}_2,\widetilde{e}_3)=-\frac{\alpha\gamma^2}{2},~~~(\nabla^k_{\widetilde{e}_2}\widetilde{\rho}^k)(\widetilde{e}_1,\widetilde{e}_3)=-\frac{\alpha\gamma^2}{2},\nonumber\\
&(\nabla^k_{\widetilde{e}_1}\widetilde{\rho}^k)(\widetilde{e}_3,\widetilde{e}_1)=0,~~~(\nabla^k_{\widetilde{e}_3}\widetilde{\rho}^k)(\widetilde{e}_1,\widetilde{e}_1)=0,~~~(\nabla^k_{\widetilde{e}_1}\widetilde{\rho}^k)(\widetilde{e}_3,\widetilde{e}_2)=-\frac{\alpha\gamma^2}{2},\nonumber\\
&(\nabla^k_{\widetilde{e}_3}\widetilde{\rho}^k)(\widetilde{e}_1,\widetilde{e}_2)=\gamma^2(\beta-\alpha),~~~(\nabla^k_{\widetilde{e}_1}\widetilde{\rho}^k)(\widetilde{e}_3,\widetilde{e}_3)=0,~~~(\nabla^k_{\widetilde{e}_3}\widetilde{\rho}^k)(\widetilde{e}_1,\widetilde{e}_3)=\frac{\alpha\beta\gamma}{2},\nonumber\\
&(\nabla^k_{\widetilde{e}_2}\widetilde{\rho}^k)(\widetilde{e}_3,\widetilde{e}_1)=-\frac{\alpha\gamma^2}{2},~~~(\nabla^k_{\widetilde{e}_3}\widetilde{\rho}^k)(\widetilde{e}_2,\widetilde{e}_1)=\gamma^2(\beta-\alpha),~~~(\nabla^k_{\widetilde{e}_2}\widetilde{\rho}^k)(\widetilde{e}_3,\widetilde{e}_2)=0,\nonumber\\
&(\nabla^k_{\widetilde{e}_3}\widetilde{\rho}^k)(\widetilde{e}_2,\widetilde{e}_2)=0,~~~(\nabla^k_{\widetilde{e}_2}\widetilde{\rho}^k)(\widetilde{e}_3,\widetilde{e}_3)=0,~~~(\nabla^k_{\widetilde{e}_3}\widetilde{\rho}^k)(\widetilde{e}_2,\widetilde{e}_3)=0.
\end{align}
Then, if $\widetilde{\rho}^k$  is a Codazzi tensor on $(G_2,\nabla^k)$, by (2.6) and (2.7), we have the following three equations:
\begin{eqnarray}
       \begin{cases}
      \beta\gamma(\alpha-\beta)=0\\[2pt]
            \gamma^2(\frac{\alpha}{2}-\beta)=0\\[2pt]
       \frac{\alpha\beta\gamma}{2}=0\\[2pt]
       \end{cases}
\end{eqnarray}
By solving (4.18), we get\\
\begin{thm}
$\widetilde{\rho}^k$  is a Codazzi tensor on $(G_2,\nabla^k)$ if and only if $\gamma\neq0,~~~\alpha=\beta=0$.
\end{thm}
\vskip 0.5 true cm
\noindent{\bf 4.3 Codazzi tensors of $G_3$}\\
\vskip 0.5 true cm
\begin{lem} (\cite{wy}) The canonical connection $\nabla^c$ of $G_3$ is given by
\begin{align}
&\nabla^c_{\widetilde{e}_1}\widetilde{e}_1=0,~~~\nabla^c_{\widetilde{e}_1}\widetilde{e}_2=0,~~~\nabla^c_{\widetilde{e}_1}\widetilde{e}_3=0,\nonumber\\
&\nabla^c_{\widetilde{e}_2}\widetilde{e}_1=0,~~~\nabla^c_{\widetilde{e}_2}\widetilde{e}_2=0,~~~\nabla^c_{\widetilde{e}_2}\widetilde{e}_3=0,\nonumber\\
&\nabla^c_{\widetilde{e}_3}\widetilde{e}_1=m_3\widetilde{e}_2,~~~\nabla^c_{\widetilde{e}_3}\widetilde{e}_2=-m_3\widetilde{e}_1,~~~\nabla^c_{\widetilde{e}_3}\widetilde{e}_3=0,
\end{align}
where
\begin{align}
m_1=\frac{\alpha-\beta-\gamma}{2},~~~m_2=\frac{\alpha-\beta+\gamma}{2},~~~m_3=\frac{\alpha+\beta-\gamma}{2}.
\end{align}
\end{lem}
Then,
\begin{align}
&\widetilde{\rho}^c(\widetilde{e}_1,\widetilde{e}_1)=-m_3\gamma,~~~\widetilde{\rho}^c(\widetilde{e}_1,\widetilde{e}_2)=\widetilde{\rho}^c(\widetilde{e}_1,\widetilde{e}_3)=0,\nonumber\\
&\widetilde{\rho}^c(\widetilde{e}_2,\widetilde{e}_2)=-m_3\gamma,~~~\widetilde{\rho}^c(\widetilde{e}_2,\widetilde{e}_3)=\widetilde{\rho}^c(\widetilde{e}_3,\widetilde{e}_3)=0.
\end{align}
By (2.5), we have
\begin{align}
&(\nabla^c_{\widetilde{e}_1}\widetilde{\rho}^c)(\widetilde{e}_2,\widetilde{e}_j)=(\nabla^c_{\widetilde{e}_2}\widetilde{\rho}^c)(\widetilde{e}_1,\widetilde{e}_j)=(\nabla^c_{\widetilde{e}_1}\widetilde{\rho}^c)(\widetilde{e}_3,\widetilde{e}_j)=0,\nonumber\\
&(\nabla^c_{\widetilde{e}_3}\widetilde{\rho}^c)(\widetilde{e}_1,\widetilde{e}_j)=(\nabla^c_{\widetilde{e}_2}\widetilde{\rho}^c)(\widetilde{e}_3,\widetilde{e}_j)=(\nabla^c_{\widetilde{e}_3}\widetilde{\rho}^c)(\widetilde{e}_2,\widetilde{e}_j)=0,
\end{align}
where $1\leq j\leq3$.\\
Then, we get
\begin{thm}
$\widetilde{\rho}^c$ is a Codazzi tensor on $(G_3,\nabla^c)$.
\end{thm}
\begin{lem} (\cite{wy}) The Kobayashi-Nomizu connection $\nabla^k$ of $G_3$ is given by
\begin{align}
&\nabla^k_{\widetilde{e}_1}\widetilde{e}_1=0,~~~\nabla^k_{\widetilde{e}_1}\widetilde{e}_2=0,~~~\nabla^k_{\widetilde{e}_1}\widetilde{e}_3=0,\nonumber\\
&\nabla^k_{\widetilde{e}_2}\widetilde{e}_1=0,~~~\nabla^k_{\widetilde{e}_2}\widetilde{e}_2=0,~~~\nabla^k_{\widetilde{e}_2}\widetilde{e}_3=0,\nonumber\\
&\nabla^k_{\widetilde{e}_3}\widetilde{e}_1=(m_3-m_1)\widetilde{e}_2,~~~\nabla^k_{\widetilde{e}_3}\widetilde{e}_2=-(m_2+m_3)\widetilde{e}_1,~~~\nabla^k_{\widetilde{e}_3}\widetilde{e}_3=0.
\end{align}
where
\begin{align}
m_1=\frac{\alpha-\beta-\gamma}{2},~~~m_2=\frac{\alpha-\beta+\gamma}{2},~~~m_3=\frac{\alpha+\beta-\gamma}{2}.
\end{align}
\end{lem}
Then,
\begin{align}
&\widetilde{\rho}^k(\widetilde{e}_1,\widetilde{e}_1)=\gamma(m_1-m_3),~~~\widetilde{\rho}^k(\widetilde{e}_1,\widetilde{e}_2)=\widetilde{\rho}^k(\widetilde{e}_1,\widetilde{e}_3)=0,\nonumber\\
&\widetilde{\rho}^k(\widetilde{e}_2,\widetilde{e}_2)=-\gamma(m_2+m_3),~~~\widetilde{\rho}^k(\widetilde{e}_2,\widetilde{e}_3)=\widetilde{\rho}^k(\widetilde{e}_3,\widetilde{e}_3)=0.
\end{align}
By (2.5), we have
\begin{align}
&(\nabla^k_{\widetilde{e}_1}\widetilde{\rho}^k)(\widetilde{e}_2,\widetilde{e}_j)=(\nabla^k_{\widetilde{e}_2}\widetilde{\rho}^k)(\widetilde{e}_1,\widetilde{e}_j)=(\nabla^k_{\widetilde{e}_1}\widetilde{\rho}^k)(\widetilde{e}_3,\widetilde{e}_j)=0,\nonumber\\
&(\nabla^k_{\widetilde{e}_3}\widetilde{\rho}^k)(\widetilde{e}_1,\widetilde{e}_j)=(\nabla^k_{\widetilde{e}_2}\widetilde{\rho}^k)(\widetilde{e}_3,\widetilde{e}_j)=(\nabla^k_{\widetilde{e}_3}\widetilde{\rho}^k)(\widetilde{e}_2,\widetilde{e}_j)=0,
\end{align}
where $1\leq j\leq3$.\\
Then, we get
\begin{thm}
$\widetilde{\rho}^k$ is a Codazzi tensor on $(G_3,\nabla^k)$.
\end{thm}
\vskip 0.5 true cm
\noindent{\bf 4.4 Codazzi tensors of $G_4$}\\
\vskip 0.5 true cm
\begin{lem} (\cite{wy}) The canonical connection $\nabla^c$ of $G_4$ is given by
\begin{align}
&\nabla^c_{\widetilde{e}_1}\widetilde{e}_1=0,~~~\nabla^c_{\widetilde{e}_1}\widetilde{e}_2=0,~~~\nabla^c_{\widetilde{e}_1}\widetilde{e}_3=0,\nonumber\\
&\nabla^c_{\widetilde{e}_2}\widetilde{e}_1=\widetilde{e}_2,~~~\nabla^c_{\widetilde{e}_2}\widetilde{e}_2=-\widetilde{e}_1,~~~\nabla^c_{\widetilde{e}_2}\widetilde{e}_3=0,\nonumber\\
&\nabla^c_{\widetilde{e}_3}\widetilde{e}_1=n_3\widetilde{e}_2,~~~\nabla^c_{\widetilde{e}_3}\widetilde{e}_2=-n_3\widetilde{e}_1,~~~\nabla^c_{\widetilde{e}_3}\widetilde{e}_3=0.
\end{align}
where
\begin{align}
n_1=\frac{\alpha}{2}+\eta-\beta,~~~n_2=\frac{\alpha}{2}-\eta,~~~n_3=\frac{\alpha}{2}+\eta.
\end{align}
\end{lem}
Then,
\begin{align}
&\widetilde{\rho}^c(\widetilde{e}_1,\widetilde{e}_1)=(2\eta-\beta)n_3-1,~~~\widetilde{\rho}^c(\widetilde{e}_1,\widetilde{e}_2)=0,~~~\widetilde{\rho}^c(\widetilde{e}_1,\widetilde{e}_3)=0,\nonumber\\
&\widetilde{\rho}^c(\widetilde{e}_2,\widetilde{e}_2)=(2\eta-\beta)n_3-1,~~~\widetilde{\rho}^c(\widetilde{e}_2,\widetilde{e}_3)=(\frac{n_3-\beta}{2},~~~\widetilde{\rho}^c(\widetilde{e}_3,\widetilde{e}_3)=0.
\end{align}
By (2.5), we have
\begin{align}
&(\nabla^c_{\widetilde{e}_1}\widetilde{\rho}^c)(\widetilde{e}_2,\widetilde{e}_1)=(\nabla^c_{\widetilde{e}_2}\widetilde{\rho}^c)(\widetilde{e}_1,\widetilde{e}_1)=(\nabla^c_{\widetilde{e}_1}\widetilde{\rho}^c)(\widetilde{e}_2,\widetilde{e}_2)=0,\nonumber\\
&(\nabla^c_{\widetilde{e}_2}\widetilde{\rho}^c)(\widetilde{e}_1,\widetilde{e}_2)=(\nabla^c_{\widetilde{e}_1}\widetilde{\rho}^c)(\widetilde{e}_2,\widetilde{e}_3)=(\nabla^c_{\widetilde{e}_2}\widetilde{\rho}^c)(\widetilde{e}_1,\widetilde{e}_3)=\frac{\beta-n_3}{2},\nonumber\\
&(\nabla^c_{\widetilde{e}_1}\widetilde{\rho}^c)(\widetilde{e}_3,\widetilde{e}_1)=(\nabla^c_{\widetilde{e}_3}\widetilde{\rho}^c)(\widetilde{e}_1,\widetilde{e}_1)=(\nabla^c_{\widetilde{e}_1}\widetilde{\rho}^c)(\widetilde{e}_3,\widetilde{e}_2)=0,\nonumber\\
&(\nabla^c_{\widetilde{e}_3}\widetilde{\rho}^c)(\widetilde{e}_1,\widetilde{e}_2)=(\nabla^c_{\widetilde{e}_1}\widetilde{\rho}^c)(\widetilde{e}_3,\widetilde{e}_3)=0,~~~(\nabla^c_{\widetilde{e}_3}\widetilde{\rho}^c)(\widetilde{e}_1,\widetilde{e}_3)=\frac{n_3(n_3-\beta)}{2},\nonumber\\
&(\nabla^c_{\widetilde{e}_2}\widetilde{\rho}^c)(\widetilde{e}_3,\widetilde{e}_1)=\frac{\beta-n_3}{2},~~~(\nabla^c_{\widetilde{e}_3}\widetilde{\rho}^c)(\widetilde{e}_2,\widetilde{e}_1)=(\nabla^c_{\widetilde{e}_2}\widetilde{\rho}^c)(\widetilde{e}_3,\widetilde{e}_2)=0,\nonumber\\
&(\nabla^c_{\widetilde{e}_3}\widetilde{\rho}^c)(\widetilde{e}_2,\widetilde{e}_2)=0,~~~(\nabla^c_{\widetilde{e}_2}\widetilde{\rho}^c)(\widetilde{e}_3,\widetilde{e}_3)=(\nabla^c_{\widetilde{e}_3}\widetilde{\rho}^c)(\widetilde{e}_2,\widetilde{e}_3)=0.
\end{align}
Then, if $\widetilde{\rho}^c$  is a Codazzi tensor on $(G_4,\nabla^c)$ , by (2.6) and (2.7), we have the following two equations:
\begin{eqnarray}
       \begin{cases}
     \frac{\beta-n_3}{2}=0\\[2pt]
      \frac{n_3(\beta-n_3)}{2}=0\\[2pt]
       \end{cases}
\end{eqnarray}
By solving (4.31), we get\\
\begin{thm}
$\widetilde{\rho}^c$  is a Codazzi tensor on $(G_4,\nabla^c)$ if and only if $\frac{\alpha}{2}+\eta-\beta=0$.
\end{thm}
\begin{lem} (\cite{wy}) The Kobayashi-Nomizu connection $\nabla^k$ of $G_4$ is given by
\begin{align}
&\nabla^k_{\widetilde{e}_1}\widetilde{e}_1=0,~~~\nabla^k_{\widetilde{e}_1}\widetilde{e}_2=0,~~~\nabla^k_{\widetilde{e}_1}\widetilde{e}_3=\widetilde{e}_3,\nonumber\\
&\nabla^k_{\widetilde{e}_2}\widetilde{e}_1=\widetilde{e}_2,~~~\nabla^k_{\widetilde{e}_2}\widetilde{e}_2=-\widetilde{e}_1,~~~\nabla^k_{\widetilde{e}_2}\widetilde{e}_3=0,\nonumber\\
&\nabla^k_{\widetilde{e}_3}\widetilde{e}_1=(n_3-n_1)\widetilde{e}_2,~~~\nabla^k_{\widetilde{e}_3}\widetilde{e}_2=-(n_2+n_3)\widetilde{e}_1,~~~\nabla^k_{\widetilde{e}_3}\widetilde{e}_3=0.
\end{align}
where
\begin{align}
n_1=\frac{\alpha}{2}+\eta-\beta,~~~n_2=\frac{\alpha}{2}-\eta,~~~n_3=\frac{\alpha}{2}+\eta.
\end{align}
\end{lem}
Then,
\begin{align}
&\widetilde{\rho}^k(\widetilde{e}_1,\widetilde{e}_1)=-[1+(\beta-2\eta)(n_3-n_1)],~~~\widetilde{\rho}^k(\widetilde{e}_1,\widetilde{e}_2)=0,~~~\widetilde{\rho}^k(\widetilde{e}_1,\widetilde{e}_3)=0,\nonumber\\
&\widetilde{\rho}^k(\widetilde{e}_2,\widetilde{e}_2)=-[1+(\beta-2\eta)(n_2+n_3)],~~~\widetilde{\rho}^k(\widetilde{e}_2,\widetilde{e}_3)=\frac{\alpha+n_3-n_1-\beta}{2},~~~\widetilde{\rho}^k(\widetilde{e}_3,\widetilde{e}_3)=0.
\end{align}
By (2.5), we have
\begin{align}
&(\nabla^k_{\widetilde{e}_1}\widetilde{\rho}^k)(\widetilde{e}_2,\widetilde{e}_1)=(\nabla^k_{\widetilde{e}_2}\widetilde{\rho}^k)(\widetilde{e}_1,\widetilde{e}_1)=(\nabla^k_{\widetilde{e}_1}\widetilde{\rho}^k)(\widetilde{e}_2,\widetilde{e}_2)=0,~~~(\nabla^k_{\widetilde{e}_2}\widetilde{\rho}^k)(\widetilde{e}_1,\widetilde{e}_2)=(n_1+n_2)(\beta-2\eta),\nonumber\\
&(\nabla^k_{\widetilde{e}_1}\widetilde{\rho}^k)(\widetilde{e}_2,\widetilde{e}_3)=\frac{n_1+\beta-\alpha-n_3}{2},~~~(\nabla^k_{\widetilde{e}_2}\widetilde{\rho}^k)(\widetilde{e}_1,\widetilde{e}_3)=\frac{n_1+\beta-\alpha-n_3}{2},\nonumber\\
&(\nabla^k_{\widetilde{e}_1}\widetilde{\rho}^k)(\widetilde{e}_3,\widetilde{e}_1)=(\nabla^k_{\widetilde{e}_3}\widetilde{\rho}^k)(\widetilde{e}_1,\widetilde{e}_1)=0,~~~(\nabla^k_{\widetilde{e}_1}\widetilde{\rho}^k)(\widetilde{e}_3,\widetilde{e}_2)=\frac{n_1+\beta-\alpha-n_3}{2},\nonumber\\
&(\nabla^k_{\widetilde{e}_3}\widetilde{\rho}^k)(\widetilde{e}_1,\widetilde{e}_2)=-(n_1+n_2),~~~(\nabla^k_{\widetilde{e}_1}\widetilde{\rho}^k)(\widetilde{e}_3,\widetilde{e}_3)=0,~~~(\nabla^k_{\widetilde{e}_3}\widetilde{\rho}^k)(\widetilde{e}_1,\widetilde{e}_3)=(n_3-n_1)\frac{n_1+\beta-\alpha-n_3}{2},\nonumber\\
&(\nabla^k_{\widetilde{e}_2}\widetilde{\rho}^k)(\widetilde{e}_3,\widetilde{e}_1)=\frac{n_1+\beta-\alpha-n_3}{2},~~~(\nabla^k_{\widetilde{e}_3}\widetilde{\rho}^k)(\widetilde{e}_2,\widetilde{e}_1)=-(n_1+n_2),~~~(\nabla^k_{\widetilde{e}_2}\widetilde{\rho}^k)(\widetilde{e}_3,\widetilde{e}_2)=0,\nonumber\\
&(\nabla^k_{\widetilde{e}_3}\widetilde{\rho}^k)(\widetilde{e}_2,\widetilde{e}_2)=(\nabla^k_{\widetilde{e}_2}\widetilde{\rho}^k)(\widetilde{e}_3,\widetilde{e}_3)=(\nabla^k_{\widetilde{e}_3}\widetilde{\rho}^k)(\widetilde{e}_2,\widetilde{e}_3)=0.
\end{align}
Then, if $\widetilde{\rho}^k$  is a Codazzi tensor on $(G_4,\nabla^k)$, by (2.6) and (2.7), we have the following three equations:
\begin{eqnarray}
       \begin{cases}
      (2\eta-\beta)(n_1+n_2)=0\\[2pt]
           \frac{3n_1+\beta-\alpha-n_3}{2}+n_2=0\\[2pt]
       (n_3-n_1)\frac{\alpha+n_3-n_1-\beta}{2}=0\\[2pt]
       \end{cases}
\end{eqnarray}
By solving (4.36) , we get\\
\begin{thm}
$\widetilde{\rho}^k$  is a Codazzi tensor on $(G_4,\nabla^k)$ if and only if $\alpha=\beta=0$.
\end{thm}
\vskip 0.5 true cm
\noindent{\bf 4.5 Codazzi tensors of $G_5$}\\
\vskip 0.5 true cm
\begin{lem} (\cite{wy}) The canonical connection $\nabla^c$ of $G_5$ is given by
\begin{align}
&\nabla^c_{\widetilde{e}_1}\widetilde{e}_1=0,~~~\nabla^c_{\widetilde{e}_1}\widetilde{e}_2=0,~~~\nabla^c_{\widetilde{e}_1}\widetilde{e}_3=0,\nonumber\\
&\nabla^c_{\widetilde{e}_2}\widetilde{e}_1=0,~~~\nabla^c_{\widetilde{e}_2}\widetilde{e}_2=0,~~~\nabla^c_{\widetilde{e}_2}\widetilde{e}_3=0,\nonumber\\
&\nabla^c_{\widetilde{e}_3}\widetilde{e}_1=\frac{\gamma-\beta}{2}\widetilde{e}_2,~~~\nabla^c_{\widetilde{e}_3}\widetilde{e}_2=\frac{\beta-\gamma}{2}\widetilde{e}_1,~~~\nabla^c_{\widetilde{e}_3}\widetilde{e}_3=0,
\end{align}
\end{lem}
Then,
\begin{align}
&\widetilde{\rho}^c(\widetilde{e}_1,\widetilde{e}_1)=\widetilde{\rho}^c(\widetilde{e}_1,\widetilde{e}_2)=\widetilde{\rho}^c(\widetilde{e}_1,\widetilde{e}_3)=0,\nonumber\\
&\widetilde{\rho}^c(\widetilde{e}_2,\widetilde{e}_2)=\widetilde{\rho}^c(\widetilde{e}_2,\widetilde{e}_3)=\widetilde{\rho}^c(\widetilde{e}_3,\widetilde{e}_3)=0.
\end{align}
By (2.5), we have
\begin{align}
&(\nabla^c_{\widetilde{e}_1}\widetilde{\rho}^c)(\widetilde{e}_2,\widetilde{e}_j)=(\nabla^c_{\widetilde{e}_2}\widetilde{\rho}^c)(\widetilde{e}_1,\widetilde{e}_j)=(\nabla^c_{\widetilde{e}_1}\widetilde{\rho}^c)(\widetilde{e}_3,\widetilde{e}_j)=0,\nonumber\\
&(\nabla^c_{\widetilde{e}_3}\widetilde{\rho}^c)(\widetilde{e}_1,\widetilde{e}_j)=(\nabla^c_{\widetilde{e}_2}\widetilde{\rho}^c)(\widetilde{e}_3,\widetilde{e}_j)=(\nabla^c_{\widetilde{e}_3}\widetilde{\rho}^c)(\widetilde{e}_2,\widetilde{e}_j)=0,
\end{align}
where $1\leq j\leq3$.\\
Then, we get
\begin{thm}
$\widetilde{\rho}^c$ is a Codazzi tensor on $(G_5,\nabla^c)$.
\end{thm}
\begin{lem} (\cite{wy}) The Kobayashi-Nomizu connection $\nabla^k$ of $G_5$ is given by
\begin{align}
&\nabla^k_{\widetilde{e}_1}\widetilde{e}_1=0,~~~\nabla^k_{\widetilde{e}_1}\widetilde{e}_2=0,~~~\nabla^k_{\widetilde{e}_1}\widetilde{e}_3=0,\nonumber\\
&\nabla^k_{\widetilde{e}_2}\widetilde{e}_1=0,~~~\nabla^k_{\widetilde{e}_2}\widetilde{e}_2=0,~~~\nabla^k_{\widetilde{e}_2}\widetilde{e}_3=0,\nonumber\\
&\nabla^k_{\widetilde{e}_3}\widetilde{e}_1=-\alpha\widetilde{e}_1-\beta\widetilde{e}_2,~~~\nabla^k_{\widetilde{e}_3}\widetilde{e}_2=-\gamma\widetilde{e}_1-\delta\widetilde{e}_2,~~~\nabla^k_{\widetilde{e}_3}\widetilde{e}_3=0.
\end{align}
\end{lem}
Then,
\begin{align}
&\widetilde{\rho}^k(\widetilde{e}_1,\widetilde{e}_1)=\widetilde{\rho}^k(\widetilde{e}_1,\widetilde{e}_2)=\widetilde{\rho}^k(\widetilde{e}_1,\widetilde{e}_3)=0,\nonumber\\
&\widetilde{\rho}^k(\widetilde{e}_2,\widetilde{e}_2)=\widetilde{\rho}^k(\widetilde{e}_2,\widetilde{e}_3)=\widetilde{\rho}^k(\widetilde{e}_3,\widetilde{e}_3)=0.
\end{align}
Then, we get
\begin{thm}
$\widetilde{\rho}^k$ is a Codazzi tensor on $(G_5,\nabla^k)$.
\end{thm}
\vskip 0.5 true cm
\noindent{\bf 4.6 Codazzi tensors of $G_6$}\\
\vskip 0.5 true cm
\begin{lem} (\cite{wy}) The canonical connection $\nabla^c$ of $G_6$ is given by
\begin{align}
&\nabla^c_{\widetilde{e}_1}\widetilde{e}_1=0,~~~\nabla^c_{\widetilde{e}_1}\widetilde{e}_2=\alpha\widetilde{e}_1,~~~\nabla^c_{\widetilde{e}_1}\widetilde{e}_3=0,\nonumber\\
&\nabla^c_{\widetilde{e}_2}\widetilde{e}_1=-\alpha\widetilde{e}_2,~~~\nabla^c_{\widetilde{e}_2}\widetilde{e}_2=\alpha\widetilde{e}_1,~~~\nabla^c_{\widetilde{e}_2}\widetilde{e}_3=0,\nonumber\\
&\nabla^c_{\widetilde{e}_3}\widetilde{e}_1=\frac{\beta-\gamma}{2}\widetilde{e}_2,~~~\nabla^c_{\widetilde{e}_3}\widetilde{e}_2=-\frac{\beta-\gamma}{2}\widetilde{e}_1,~~~\nabla^c_{\widetilde{e}_3}\widetilde{e}_3=0.
\end{align}
\end{lem}
Then,
\begin{align}
&\widetilde{\rho}^c(\widetilde{e}_1,\widetilde{e}_1)=\frac{1}{2}\beta(\beta-\gamma)-\alpha^2,~~~\widetilde{\rho}^c(\widetilde{e}_1,\widetilde{e}_2)=0,~~~\widetilde{\rho}^c(\widetilde{e}_1,\widetilde{e}_3)=0,\nonumber\\
&\widetilde{\rho}^c(\widetilde{e}_2,\widetilde{e}_2)=\frac{1}{2}\beta(\beta-\gamma)-\alpha^2,~~~\widetilde{\rho}^c(\widetilde{e}_2,\widetilde{e}_3)=\frac{1}{2}[-\alpha\gamma+\frac{1}{2}\delta(\beta-\gamma)],~~~\widetilde{\rho}^c(\widetilde{e}_3,\widetilde{e}_3)=0.
\end{align}
By (2.5), we have
\begin{align}
&(\nabla^c_{\widetilde{e}_1}\widetilde{\rho}^c)(\widetilde{e}_2,\widetilde{e}_1)=0,~~~(\nabla^c_{\widetilde{e}_2}\widetilde{\rho}^c)(\widetilde{e}_1,\widetilde{e}_1)=0,~~~(\nabla^c_{\widetilde{e}_1}\widetilde{\rho}^c)(\widetilde{e}_2,\widetilde{e}_2)=0,\nonumber\\
&(\nabla^c_{\widetilde{e}_2}\widetilde{\rho}^c)(\widetilde{e}_1,\widetilde{e}_2)=0,~~~(\nabla^c_{\widetilde{e}_1}\widetilde{\rho}^c)(\widetilde{e}_2,\widetilde{e}_3)=0,~~~(\nabla^c_{\widetilde{e}_2}\widetilde{\rho}^c)(\widetilde{e}_1,\widetilde{e}_3)=\frac{\alpha}{2}[-\alpha\gamma+\frac{1}{2}\delta(\beta-\gamma)],\nonumber\\
&(\nabla^c_{\widetilde{e}_1}\widetilde{\rho}^c)(\widetilde{e}_3,\widetilde{e}_1)=0,~~~(\nabla^c_{\widetilde{e}_3}\widetilde{\rho}^c)(\widetilde{e}_1,\widetilde{e}_1)=0,~~~(\nabla^c_{\widetilde{e}_1}\widetilde{\rho}^c)(\widetilde{e}_3,\widetilde{e}_2)=0,\nonumber\\
&(\nabla^c_{\widetilde{e}_3}\widetilde{\rho}^c)(\widetilde{e}_1,\widetilde{e}_2)=0,~~~(\nabla^c_{\widetilde{e}_1}\widetilde{\rho}^c)(\widetilde{e}_3,\widetilde{e}_3)=0,~~~(\nabla^c_{\widetilde{e}_3}\widetilde{\rho}^c)(\widetilde{e}_1,\widetilde{e}_3)=\frac{\gamma-\beta}{4}[-\alpha\gamma+\frac{1}{2}\delta(\beta-\gamma)],\nonumber\\
&(\nabla^c_{\widetilde{e}_2}\widetilde{\rho}^c)(\widetilde{e}_3,\widetilde{e}_1)=\frac{\alpha}{2}[-\alpha\gamma+\frac{1}{2}\delta(\beta-\gamma)],~~~(\nabla^c_{\widetilde{e}_3}\widetilde{\rho}^c)(\widetilde{e}_2,\widetilde{e}_1)=0,~~~(\nabla^c_{\widetilde{e}_2}\widetilde{\rho}^c)(\widetilde{e}_3,\widetilde{e}_2)=0,\nonumber\\
&(\nabla^c_{\widetilde{e}_3}\widetilde{\rho}^c)(\widetilde{e}_2,\widetilde{e}_2)=0,~~~(\nabla^c_{\widetilde{e}_2}\widetilde{\rho}^c)(\widetilde{e}_3,\widetilde{e}_3)=0,~~~(\nabla^c_{\widetilde{e}_3}\widetilde{\rho}^c)(\widetilde{e}_2,\widetilde{e}_3)=0.
\end{align}
Then, if $\widetilde{\rho}^c$  is a Codazzi tensor on $(G_6,\nabla^c)$, by (2.6) and (2.7), we have the following two equations:
\begin{eqnarray}
       \begin{cases}
      \frac{\alpha}{2}[-\alpha\gamma+\frac{1}{2}\delta(\beta-\gamma)]=0\\[2pt]
      \frac{\beta-\gamma}{4}[-\alpha\gamma+\frac{1}{2}\delta(\beta-\gamma)]=0\\[2pt]
       \end{cases}
\end{eqnarray}
By solving (4.45), we get\\
\begin{thm}
$\widetilde{\rho}^c$  is a Codazzi tensor on $(G_6,\nabla^c)$ if and only if
\begin{eqnarray*}
&&(1)\alpha=\beta=\gamma=0,~~~\delta\neq0;\\
&&(2)\alpha\neq0,~~~\beta=\gamma=0,~~~\alpha+\delta\neq0;\\
&&(3)\alpha\neq0,~~~\delta=\gamma=0.
\end{eqnarray*}
\end{thm}
\begin{lem} (\cite{wy}) The Kobayashi-Nomizu connection $\nabla^k$ of $G_6$ is given by
\begin{align}
&\nabla^k_{\widetilde{e}_1}\widetilde{e}_1=0,~~~\nabla^k_{\widetilde{e}_1}\widetilde{e}_2=0,~~~\nabla^k_{\widetilde{e}_1}\widetilde{e}_3=\delta\widetilde{e}_3,\nonumber\\
&\nabla^k_{\widetilde{e}_2}\widetilde{e}_1=-\alpha\widetilde{e}_2,~~~\nabla^k_{\widetilde{e}_2}\widetilde{e}_2=\alpha\widetilde{e}_1,~~~\nabla^k_{\widetilde{e}_2}\widetilde{e}_3=0,\nonumber\\
&\nabla^k_{\widetilde{e}_3}\widetilde{e}_1=-\gamma\widetilde{e}_2,~~~\nabla^k_{\widetilde{e}_3}\widetilde{e}_2=0,~~~\nabla^k_{\widetilde{e}_3}\widetilde{e}_3=0.
\end{align}
\end{lem}
Then,
\begin{align}
&\widetilde{\rho}^k(\widetilde{e}_1,\widetilde{e}_1)=-(\alpha^2+\beta\gamma),~~~\widetilde{\rho}^k(\widetilde{e}_1,\widetilde{e}_2)=0,~~~\widetilde{\rho}^k(\widetilde{e}_1,\widetilde{e}_3)=0,\nonumber\\
&\widetilde{\rho}^k(\widetilde{e}_2,\widetilde{e}_2)=-\alpha^2,~~~\widetilde{\rho}^k(\widetilde{e}_2,\widetilde{e}_3)=-\frac{\alpha\gamma}{2},~~~\widetilde{\rho}^k(\widetilde{e}_3,\widetilde{e}_3)=0.
\end{align}
By (2.5), we have
\begin{align}
&(\nabla^k_{\widetilde{e}_1}\widetilde{\rho}^k)(\widetilde{e}_2,\widetilde{e}_1)=(\nabla^k_{\widetilde{e}_2}\widetilde{\rho}^k)(\widetilde{e}_1,\widetilde{e}_1)=(\nabla^k_{\widetilde{e}_1}\widetilde{\rho}^k)(\widetilde{e}_2,\widetilde{e}_2)=0,\nonumber\\
&(\nabla^k_{\widetilde{e}_2}\widetilde{\rho}^k)(\widetilde{e}_1,\widetilde{e}_2)=\alpha\beta\gamma,~~~(\nabla^k_{\widetilde{e}_1}\widetilde{\rho}^k)(\widetilde{e}_2,\widetilde{e}_3)=(\nabla^k_{\widetilde{e}_2}\widetilde{\rho}^k)(\widetilde{e}_1,\widetilde{e}_3)=0,\nonumber\\
&(\nabla^k_{\widetilde{e}_1}\widetilde{\rho}^k)(\widetilde{e}_3,\widetilde{e}_1)=(\nabla^k_{\widetilde{e}_3}\widetilde{\rho}^k)(\widetilde{e}_1,\widetilde{e}_1)=0,~~~(\nabla^k_{\widetilde{e}_1}\widetilde{\rho}^k)(\widetilde{e}_3,\widetilde{e}_2)=-\frac{\alpha\gamma^2}{2},\nonumber\\
&(\nabla^k_{\widetilde{e}_3}\widetilde{\rho}^k)(\widetilde{e}_1,\widetilde{e}_2)=-\alpha^2\gamma,~~~(\nabla^k_{\widetilde{e}_1}\widetilde{\rho}^k)(\widetilde{e}_3,\widetilde{e}_3)=(\nabla^k_{\widetilde{e}_3}\widetilde{\rho}^k)(\widetilde{e}_1,\widetilde{e}_3)=0,\nonumber\\
&(\nabla^k_{\widetilde{e}_2}\widetilde{\rho}^k)(\widetilde{e}_3,\widetilde{e}_1)=-\alpha^2\gamma,~~~(\nabla^k_{\widetilde{e}_3}\widetilde{\rho}^k)(\widetilde{e}_2,\widetilde{e}_1)=(\nabla^k_{\widetilde{e}_2}\widetilde{\rho}^k)(\widetilde{e}_3,\widetilde{e}_2)=0,\nonumber\\
&(\nabla^k_{\widetilde{e}_3}\widetilde{\rho}^k)(\widetilde{e}_2,\widetilde{e}_2)=(\nabla^k_{\widetilde{e}_2}\widetilde{\rho}^k)(\widetilde{e}_3,\widetilde{e}_3)=(\nabla^k_{\widetilde{e}_3}\widetilde{\rho}^k)(\widetilde{e}_2,\widetilde{e}_3)=0.
\end{align}
Then, if $\widetilde{\rho}^k$  is a Codazzi tensor on $(G_6,\nabla^k)$ , by (2.6) and (2.7), we have the following two equations:
\begin{eqnarray}
       \begin{cases}
      \alpha\beta\gamma=0\\[2pt]
      \alpha^2\gamma=0\\[2pt]
       \end{cases}
\end{eqnarray}
By solving (4.49), we get\\
\begin{thm}
$\widetilde{\rho}^k$  is a Codazzi tensor on $(G_6,\nabla^k)$ if and only if
\begin{eqnarray*}
&&(1)\alpha=\beta=0,~~~\delta\neq0;\\
&&(2)\alpha\neq0,~~~\beta\delta=\gamma=0,~~~\alpha+\delta\neq0.
\end{eqnarray*}
\end{thm}
\vskip 0.5 true cm
\noindent{\bf 4.7 Codazzi tensors of $G_7$}\\
\vskip 0.5 true cm
\begin{lem} (\cite{wy}) The canonical connection $\nabla^c$ of $G_7$ is given by
\begin{align}
&\nabla^c_{\widetilde{e}_1}\widetilde{e}_1=\alpha\widetilde{e}_2,~~~\nabla^c_{\widetilde{e}_1}\widetilde{e}_2=-\alpha\widetilde{e}_1,~~~\nabla^c_{\widetilde{e}_1}\widetilde{e}_3=0,\nonumber\\
&\nabla^c_{\widetilde{e}_2}\widetilde{e}_1=\beta\widetilde{e}_2,~~~\nabla^c_{\widetilde{e}_2}\widetilde{e}_2=-\beta\widetilde{e}_1,~~~\nabla^c_{\widetilde{e}_2}\widetilde{e}_3=0,\nonumber\\
&\nabla^c_{\widetilde{e}_3}\widetilde{e}_1=(\frac{\gamma}{2}-\beta)\widetilde{e}_2,~~~\nabla^c_{\widetilde{e}_3}\widetilde{e}_2=(\beta-\frac{\gamma}{2})\widetilde{e}_1,~~~\nabla^c_{\widetilde{e}_3}\widetilde{e}_3=0.
\end{align}
\end{lem}
Then,
\begin{align}
&\widetilde{\rho}^c(\widetilde{e}_1,\widetilde{e}_1)=-(\alpha^2+\frac{\beta\gamma}{2}),~~~\widetilde{\rho}^c(\widetilde{e}_1,\widetilde{e}_2)=0,~~~\widetilde{\rho}^c(\widetilde{e}_1,\widetilde{e}_3)=-\frac{1}{2}(\alpha\gamma+\frac{\delta\gamma}{2}),\nonumber\\
&\widetilde{\rho}^c(\widetilde{e}_2,\widetilde{e}_2)=-(\alpha^2+\frac{\beta\gamma2}{2}),~~~\widetilde{\rho}^c(\widetilde{e}_2,\widetilde{e}_3)=\frac{1}{2}(\alpha^2+\frac{\beta\gamma}{2}),~~~\widetilde{\rho}^c(\widetilde{e}_3,\widetilde{e}_3)=0.
\end{align}
By (2.5), we have
\begin{align}
&(\nabla^c_{\widetilde{e}_1}\widetilde{\rho}^c)(\widetilde{e}_2,\widetilde{e}_1)=(\nabla^c_{\widetilde{e}_2}\widetilde{\rho}^c)(\widetilde{e}_1,\widetilde{e}_1)=(\nabla^c_{\widetilde{e}_1}\widetilde{\rho}^c)(\widetilde{e}_2,\widetilde{e}_2)=0,\nonumber\\
&(\nabla^c_{\widetilde{e}_2}\widetilde{\rho}^c)(\widetilde{e}_1,\widetilde{e}_2)=0,~~~(\nabla^c_{\widetilde{e}_1}\widetilde{\rho}^c)(\widetilde{e}_2,\widetilde{e}_3)=-\frac{\alpha}{2}(\alpha\gamma+\frac{\beta\gamma}{2}),~~~(\nabla^c_{\widetilde{e}_2}\widetilde{\rho}^c)(\widetilde{e}_1,\widetilde{e}_3)=-\frac{\beta}{2}(\alpha^2+\frac{\beta\gamma}{2}),\nonumber\\
&(\nabla^c_{\widetilde{e}_1}\widetilde{\rho}^c)(\widetilde{e}_3,\widetilde{e}_1)=-\frac{\alpha}{2}(\alpha^2+\frac{\beta\gamma}{2}),~~~(\nabla^c_{\widetilde{e}_3}\widetilde{\rho}^c)(\widetilde{e}_1,\widetilde{e}_1)=0,~~~(\nabla^c_{\widetilde{e}_1}\widetilde{\rho}^c)(\widetilde{e}_3,\widetilde{e}_2)=-\frac{\alpha}{2}(\alpha\gamma+\frac{\delta\gamma}{2}),\nonumber\\
&(\nabla^c_{\widetilde{e}_3}\widetilde{\rho}^c)(\widetilde{e}_1,\widetilde{e}_2)=(\nabla^c_{\widetilde{e}_1}\widetilde{\rho}^c)(\widetilde{e}_3,\widetilde{e}_3)=0,~~~(\nabla^c_{\widetilde{e}_3}\widetilde{\rho}^c)(\widetilde{e}_1,\widetilde{e}_3)=(\beta-\frac{\gamma}{2})(\frac{\alpha^2}{2}+\frac{\beta\gamma}{4}),\nonumber\\
&(\nabla^c_{\widetilde{e}_2}\widetilde{\rho}^c)(\widetilde{e}_3,\widetilde{e}_1)=-\frac{\beta}{2}(\alpha^2+\frac{\beta\gamma}{2}),~~~(\nabla^c_{\widetilde{e}_3}\widetilde{\rho}^c)(\widetilde{e}_2,\widetilde{e}_1)=0,~~~(\nabla^c_{\widetilde{e}_2}\widetilde{\rho}^c)(\widetilde{e}_3,\widetilde{e}_2)=-\frac{\beta}{2}(\alpha\gamma+\frac{\beta\gamma}{2}),\nonumber\\
&(\nabla^c_{\widetilde{e}_3}\widetilde{\rho}^c)(\widetilde{e}_2,\widetilde{e}_2)=-2\alpha^2,~~~(\nabla^c_{\widetilde{e}_2}\widetilde{\rho}^c)(\widetilde{e}_3,\widetilde{e}_3)=0,~~~(\nabla^c_{\widetilde{e}_3}\widetilde{\rho}^c)(\widetilde{e}_2,\widetilde{e}_3)=\frac{1}{2}(\beta-\frac{\gamma}{2})(\alpha\gamma+\frac{\delta\gamma}{2}).
\end{align}
Then, if $\widetilde{\rho}^c$  is a Codazzi tensor on $(G_7,\nabla^c)$, by (2.6) and (2.7), we have the following seven equations:
\begin{eqnarray}
       \begin{cases}
        \frac{\beta}{2}(\alpha^2+\frac{\beta\gamma}{2})-\frac{\alpha}{2}(\alpha\gamma+\frac{\delta\gamma}{2})=0\\[2pt]
         \frac{\alpha}{2}(\alpha^2+\frac{\beta\gamma}{2})=0\\[2pt]
          \frac{\alpha}{2}(\alpha\gamma+\frac{\delta\gamma}{2})=0\\[2pt]
           \frac{1}{2}(\frac{\gamma}{2}-\beta)(\alpha^2+\frac{\beta\gamma}{2})=0\\[2pt]
      \frac{\beta}{2}(\alpha^2+\frac{\beta\gamma}{2})=0\\[2pt]
       \frac{\beta}{2}(\alpha\gamma+\frac{\delta\gamma}{2})=0\\[2pt]
        \frac{1}{2}(\frac{\gamma}{2}-\beta)(\alpha\gamma+\frac{\delta\gamma}{2})=0\\[2pt]
       \end{cases}
\end{eqnarray}
By solving (4.53), we get\\
\begin{thm}
$\widetilde{\rho}^c$  is a Codazzi tensor on $(G_7,\nabla^c)$ if and only if $\alpha=\gamma=0,~~~\delta\neq0$.
\end{thm}
\begin{lem} (\cite{wy}) The Kobayashi-Nomizu connection $\nabla^k$ of $G_7$ is given by
\begin{align}
&\nabla^k_{\widetilde{e}_1}\widetilde{e}_1=\alpha\widetilde{e}_2,~~~\nabla^k_{\widetilde{e}_1}\widetilde{e}_2=-\alpha\widetilde{e}_1,~~~\nabla^k_{\widetilde{e}_1}\widetilde{e}_3=\beta\widetilde{e}_3,\nonumber\\
&\nabla^k_{\widetilde{e}_2}\widetilde{e}_1=\beta\widetilde{e}_2,~~~\nabla^k_{\widetilde{e}_2}\widetilde{e}_2=-\beta\widetilde{e}_1,~~~\nabla^k_{\widetilde{e}_2}\widetilde{e}_3=\delta\widetilde{e}_3,\nonumber\\
&\nabla^k_{\widetilde{e}_3}\widetilde{e}_1=-\alpha\widetilde{e}_1-\beta\widetilde{e}_2,~~~\nabla^k_{\widetilde{e}_3}\widetilde{e}_2=-\gamma\widetilde{e}_1-\delta\widetilde{e}_2,~~~\nabla^k_{\widetilde{e}_3}\widetilde{e}_3=0.
\end{align}
\end{lem}
Then,
\begin{align}
&\widetilde{\rho}^k(\widetilde{e}_1,\widetilde{e}_1)=-\alpha^2,~~~\widetilde{\rho}^k(\widetilde{e}_1,\widetilde{e}_2)=\frac{\beta}{2}(\delta-\alpha),~~~\widetilde{\rho}^k(\widetilde{e}_1,\widetilde{e}_3)=\beta(\alpha+\delta),\nonumber\\
&\widetilde{\rho}^k(\widetilde{e}_2,\widetilde{e}_2)=-(\alpha^2+\beta^2+\beta\gamma),~~~\widetilde{\rho}^k(\widetilde{e}_2,\widetilde{e}_3)=\frac{\alpha\delta+\beta\gamma+2\delta^2}{2},~~~\widetilde{\rho}^k(\widetilde{e}_3,\widetilde{e}_3)=0.
\end{align}
By (2.5), we have
\begin{align}
&(\nabla^k_{\widetilde{e}_1}\widetilde{\rho}^k)(\widetilde{e}_2,\widetilde{e}_1)=\alpha\beta(\beta+\gamma),~~~(\nabla^k_{\widetilde{e}_2}\widetilde{\rho}^k)(\widetilde{e}_1,\widetilde{e}_1)=\beta^2(\alpha-\delta),~~~(\nabla^k_{\widetilde{e}_1}\widetilde{\rho}^k)(\widetilde{e}_2,\widetilde{e}_2)=\alpha\beta(\delta-\alpha),\nonumber\\
&(\nabla^k_{\widetilde{e}_2}\widetilde{\rho}^k)(\widetilde{e}_1,\widetilde{e}_2)=\beta^2(\beta+\gamma),~~~(\nabla^k_{\widetilde{e}_1}\widetilde{\rho}^k)(\widetilde{e}_2,\widetilde{e}_3)=\alpha^2\beta+\frac{\alpha\beta\delta-\beta^2\gamma}{2}-\beta\delta^2,\nonumber\\
&(\nabla^k_{\widetilde{e}_2}\widetilde{\rho}^k)(\widetilde{e}_1,\widetilde{e}_3)=-2\beta\delta^2-\frac{\beta^2\gamma+3\alpha\beta\delta}{2},~~~(\nabla^k_{\widetilde{e}_1}\widetilde{\rho}^k)(\widetilde{e}_3,\widetilde{e}_1)=-\beta^2(\alpha+\delta)-\frac{\alpha}{2}(\alpha\delta+\beta\gamma+2\delta^2),\nonumber\\
&(\nabla^k_{\widetilde{e}_3}\widetilde{\rho}^k)(\widetilde{e}_1,\widetilde{e}_1)=\beta^2(\delta-\alpha)-2\alpha^3,~~~(\nabla^k_{\widetilde{e}_1}\widetilde{\rho}^k)(\widetilde{e}_3,\widetilde{e}_2)=-\frac{\beta}{2}(\alpha^2+\beta\gamma+2\delta^2),\nonumber\\
&(\nabla^k_{\widetilde{e}_3}\widetilde{\rho}^k)(\widetilde{e}_1,\widetilde{e}_2)=\frac{\beta(\delta^2-\alpha^2)}{2}-\beta(\alpha^2+\beta^2+\beta\gamma)-\alpha^2\gamma,~~~(\nabla^k_{\widetilde{e}_1}\widetilde{\rho}^k)(\widetilde{e}_3,\widetilde{e}_3)=0,\nonumber\\
&(\nabla^k_{\widetilde{e}_3}\widetilde{\rho}^k)(\widetilde{e}_1,\widetilde{e}_3)=\alpha\beta(\alpha+\delta)+\frac{\beta}{2}(\alpha\delta+\beta\gamma+2\delta^2),~~~(\nabla^k_{\widetilde{e}_2}\widetilde{\rho}^k)(\widetilde{e}_3,\widetilde{e}_1)=-\beta\delta(\alpha+\delta)-\frac{\beta}{2}(\alpha\delta+\beta\gamma+2\delta^2),\nonumber\\
&(\nabla^k_{\widetilde{e}_3}\widetilde{\rho}^k)(\widetilde{e}_2,\widetilde{e}_1)=\frac{\beta(\delta^2-\alpha^2)}{2}-\beta(\alpha^2+\beta^2+\beta\gamma)-\alpha^2\gamma,~~~(\nabla^k_{\widetilde{e}_2}\widetilde{\rho}^k)(\widetilde{e}_3,\widetilde{e}_2)=\beta^2(\alpha+\delta)-\frac{\delta}{2}(\alpha\delta+\beta\gamma+2\delta^2),\nonumber\\
&(\nabla^k_{\widetilde{e}_3}\widetilde{\rho}^k)(\widetilde{e}_2,\widetilde{e}_2)=\beta\gamma(\delta-\alpha)-2\delta(\alpha^2+\beta^2+\beta\gamma),~~~(\nabla^k_{\widetilde{e}_2}\widetilde{\rho}^k)(\widetilde{e}_3,\widetilde{e}_3)=0,\nonumber\\
&(\nabla^k_{\widetilde{e}_3}\widetilde{\rho}^k)(\widetilde{e}_2,\widetilde{e}_3)=\beta\gamma(\alpha+\delta)+\frac{\delta}{2}(\alpha\delta+\beta\gamma+2\delta^2).
\end{align}
Then, if $\widetilde{\rho}^k$  is a Codazzi tensor on $(G_7,\nabla^k)$, by (2.6) and (2.7), we have the following nine equations:
\begin{eqnarray}
       \begin{cases}
      \beta(\alpha\gamma+\beta\delta)=0\\[2pt]
            \beta(\alpha\delta-\alpha^2-\beta^2-\beta\gamma)=0\\[2pt]
       \beta(\alpha+\delta)^2=0\\[2pt]
             2\alpha^3-2\beta^2\delta-\alpha\delta^2-\frac{\alpha\beta\gamma+\alpha^2\delta}{2}=0\\[2pt]
      \beta(\alpha^2+\beta^2+\frac{\beta\gamma+\delta^2}{2})+\alpha^2\gamma=0\\[2pt]
      \alpha\beta(\alpha+\delta)+\frac{\beta}{2}(\alpha\delta+\beta\gamma+2\delta^2)=0\\[2pt]
            \alpha^2\gamma+\beta^3+\frac{3\alpha^2\beta+\beta^2\gamma-3\alpha\beta\delta-5\beta\delta^2}{2}=0\\[2pt]
      \alpha\beta\gamma-\delta^3+3\beta^2\delta+2\alpha^2\delta-\frac{\beta\delta\gamma-\alpha\delta^2}{2}=0\\[2pt]
      \alpha\beta\gamma+\delta^3+\frac{3\beta\delta\gamma+\alpha\delta^2}{2}=0\\[2pt]
       \end{cases}
\end{eqnarray}
By solving (4.57), we get $\alpha=\delta=0$, there is a contradiction. So\\
\begin{thm}
$\widetilde{\rho}^k$  is not a Codazzi tensor on $(G_7,\nabla^k)$.
\end{thm}
\section{Quasi-statistial structure associated to canonical connections and Kobayashi-Nomizu connections on three-dimensional Lorentzian Lie groups}
The torsion tensor of $(G_i,g,\nabla^c)$ is defined by
\begin{equation}
T^c(X,Y)=\nabla^c_XY-\nabla^c_YX-[X,Y]
\end{equation}
The torsion tensor of $(G_i,g,\nabla^k)$ is defined by
\begin{equation}
T^k(X,Y)=\nabla^k_XY-\nabla^k_YX-[X,Y]
\end{equation}
Then,
for $G_1$, we have
\begin{align}
T^c(\widetilde{e}_1,\widetilde{e}_2)=\beta \widetilde{e}_3,~~~T^c(\widetilde{e}_1,\widetilde{e}_3)=\alpha\widetilde{e}_1+\frac{\beta}{2}\widetilde{e}_2,~~~T^c(\widetilde{e}_2,\widetilde{e}_3)=-\frac{\beta}{2}\widetilde{e}_1-\alpha\widetilde{e}_2-\alpha\widetilde{e}_3.
\end{align}
\begin{align}
&\widetilde{\rho}^c(T^c(\widetilde{e}_1,\widetilde{e}_2),\widetilde{e}_1)=\frac{\alpha\beta^2}{4},~~~\widetilde{\rho}^c(T^c(\widetilde{e}_1,\widetilde{e}_2),\widetilde{e}_2)=\frac{\alpha^2\beta}{2},~~~\widetilde{\rho}^c(T^c(\widetilde{e}_1,\widetilde{e}_2),\widetilde{e}_3)=0,\nonumber\\
&\widetilde{\rho}^c(T^c(\widetilde{e}_1,\widetilde{e}_3),\widetilde{e}_1)=-\alpha(\alpha^2+\frac{\beta}{2}),~~~\widetilde{\rho}^c(T^c(\widetilde{e}_1,\widetilde{e}_3),\widetilde{e}_2)=-\frac{\beta}{2}(\alpha^2+\frac{\beta}{2}),~~~\widetilde{\rho}^c(T^c(\widetilde{e}_1,\widetilde{e}_3),\widetilde{e}_1)=\frac{\alpha^2\beta}{2},\nonumber\\
&\widetilde{\rho}^c(T^c(\widetilde{e}_2,\widetilde{e}_3),\widetilde{e}_1)=\frac{\beta}{4}(\alpha^2+\beta^2),~~~\widetilde{\rho}^c(T^c(\widetilde{e}_2,\widetilde{e}_3),\widetilde{e}_1)=\frac{\alpha}{2}(\alpha^2+\beta^2),~~~\widetilde{\rho}^c(T^c(\widetilde{e}_2,\widetilde{e}_3),\widetilde{e}_1)=-\alpha(\frac{\beta^2}{8}+\frac{\alpha^2}{2}).
\end{align}
Then, if $(G_1,\nabla^c,\widetilde{\rho}^c)$ is a quasi-statistical structure, by (3.2) and (3.3), we have the following four equations:
\begin{eqnarray}
        \begin{cases}
       \frac{\alpha\beta^2}{4}=0 \\[2pt]
              \frac{\alpha^2\beta}{2}=0 \\[2pt]
      \frac{\alpha}{2}(\alpha^2+\beta^2)=0\\[2pt]
        \frac{\beta}{4}(\alpha^2+\beta^2)=0\\[2pt]
       \end{cases}
\end{eqnarray}
By solving (5.5), we get $\alpha=0$, there is a contradiction. So\\
\begin{thm}
 $(G_1,\nabla^c,\widetilde{\rho}^c)$ is not a quasi-statistical structure.
\end{thm}
Similarly,
\begin{align}
T^k(\widetilde{e}_1,\widetilde{e}_2)=\beta \widetilde{e}_3,~~~T^k(\widetilde{e}_1,\widetilde{e}_3)=T^k(\widetilde{e}_2,\widetilde{e}_3)=0.
\end{align}
\begin{align}
&\widetilde{\rho}^k(T^k(\widetilde{e}_1,\widetilde{e}_2),\widetilde{e}_1)=0,~~~\widetilde{\rho}^k(T^k(\widetilde{e}_1,\widetilde{e}_2),\widetilde{e}_2)=\frac{\alpha^2\beta}{2},~~~\widetilde{\rho}^k(T^k(\widetilde{e}_1,\widetilde{e}_2),\widetilde{e}_3)=0,\\\nonumber
&\widetilde{\rho}^k(T^k(\widetilde{e}_1,\widetilde{e}_3),\widetilde{e}_j)=0,~~~\widetilde{\rho}^k(T^k(\widetilde{e}_2,\widetilde{e}_3),\widetilde{e}_j)=0,
\end{align}
where $1\leq j\leq 3$.\\
Then, if $(G_1,\nabla^k,\widetilde{\rho}^k)$ is a quasi-statistical structure, by (3.2) and (3.3), we have the following four equations:
\begin{eqnarray}
        \begin{cases}
       \frac{\alpha\beta^2}{2}=0 \\[2pt]
              \frac{3\alpha^2\beta}{2}=0 \\[2pt]
      \frac{3\alpha^3}{2}=0\\[2pt]
        \frac{\alpha}{2}(\alpha^2-\beta^2)=0\\[2pt]
       \end{cases}
\end{eqnarray}
By solving (5.8), we get $\alpha=0$, there is a contradiction. So\\
\begin{thm}
 $(G_1,\nabla^k,\widetilde{\rho}^k)$ is not a quasi-statistical structure.
\end{thm}
For $G_2$, we have
\begin{align}
T^c(\widetilde{e}_1,\widetilde{e}_2)=\beta \widetilde{e}_3,~~~T^c(\widetilde{e}_1,\widetilde{e}_3)=(\beta-\frac{\alpha}{2})\widetilde{e}_2+\gamma \widetilde{e}_3,~~~T^c(\widetilde{e}_2,\widetilde{e}_3)=-\frac{\alpha}{2}\widetilde{e}_1.
\end{align}
\begin{align}
&\widetilde{\rho}^c(T^c(\widetilde{e}_1,\widetilde{e}_2),\widetilde{e}_1)=0,~~~\widetilde{\rho}^c(T^c(\widetilde{e}_1,\widetilde{e}_2),\widetilde{e}_2)=\beta\gamma(\frac{\beta}{2}-\frac{\gamma}{4}),~~~\widetilde{\rho}^B(T^B(\widetilde{e}_1,\widetilde{e}_2),\widetilde{e}_3)=0,\\\nonumber
&\widetilde{\rho}^c(T^c(\widetilde{e}_1,\widetilde{e}_3),\widetilde{e}_1)=0,~~~\widetilde{\rho}^c(T^c(\widetilde{e}_1,\widetilde{e}_3),\widetilde{e}_2)=(\frac{\alpha}{4}-\frac{\beta}{2})(\alpha\beta+\gamma^2),~~~\widetilde{\rho}^c(T^c(\widetilde{e}_1,\widetilde{e}_3),\widetilde{e}_1)=(\beta-\frac{\alpha}{2})(\frac{\beta\gamma}{2}-\frac{\alpha\gamma}{4}),\\\nonumber
&\widetilde{\rho}^c(T^c(\widetilde{e}_2,\widetilde{e}_3),\widetilde{e}_1)=\frac{\alpha}{2}(\gamma^2+\frac{\alpha\beta}{2}),~~~\widetilde{\rho}^c(T^c(\widetilde{e}_2,\widetilde{e}_3),\widetilde{e}_1)=\frac{\alpha}{2}(\alpha^2+\beta^2),~~~\widetilde{\rho}^c(T^c(\widetilde{e}_2,\widetilde{e}_3),\widetilde{e}_3)=0.
\end{align}
Then, if $(G_2,\nabla^c,\widetilde{\rho}^c)$ is a quasi-statistical structure, by (3.2) and (3.3), we have the following four equations:
\begin{eqnarray}
        \begin{cases}
       \beta\gamma(\frac{\beta}{2}-\frac{\alpha}{4})=0 \\[2pt]
              \gamma^2(\frac{\beta}{2}-\frac{\alpha}{4})=0 \\[2pt]
      (\frac{\alpha}{4}-\frac{\beta}{2})(\gamma^2+\alpha\beta)=0\\[2pt]
        \frac{\alpha(\gamma^2+\alpha\beta)}{4}+\frac{\beta\gamma^2}{2}=0\\[2pt]
       \end{cases}
\end{eqnarray}
By solving (5.11), we get\\
\begin{thm}
 $(G_2,\nabla^c,\widetilde{\rho}^c)$ is a quasi-statistical structure if and only if $\gamma\neq0,~~~\alpha=\beta=0$.
\end{thm}
Similarly,
\begin{align}
T^k(\widetilde{e}_1,\widetilde{e}_2)=\beta \widetilde{e}_3,~~~T^k(\widetilde{e}_1,\widetilde{e}_3)=T^k(\widetilde{e}_2,\widetilde{e}_3)=0.
\end{align}
\begin{align}
&\widetilde{\rho}^k(T^k(\widetilde{e}_1,\widetilde{e}_2),\widetilde{e}_1)=0,~~~\widetilde{\rho}^k(T^k(\widetilde{e}_1,\widetilde{e}_2),\widetilde{e}_2)=\frac{\alpha\beta\gamma}{2},~~~\widetilde{\rho}^k(T^k(\widetilde{e}_1,\widetilde{e}_2),\widetilde{e}_3)=0,\\\nonumber
&\widetilde{\rho}^k(T^k(\widetilde{e}_1,\widetilde{e}_3),\widetilde{e}_j)=0,~~~\widetilde{\rho}^k(T^k(\widetilde{e}_2,\widetilde{e}_3),\widetilde{e}_j)=0.
\end{align}
where $1\leq j\leq3$.\\
Then, if $(G_2,\nabla^k,\widetilde{\rho}^k)$ is a quasi-statistical structure, by (3.2) and (3.3), we have the following three equations:
\begin{eqnarray}
        \begin{cases}
       \beta\gamma(\frac{\alpha}{2}-\beta)=0 \\[2pt]
              \gamma^2(\frac{\alpha}{2}-\beta)=0 \\[2pt]
      \frac{\alpha\beta\gamma}{2}=0\\[2pt]
       \end{cases}
\end{eqnarray}
By solving (5.14), we get\\
\begin{thm}
 $(G_2,\nabla^k,\widetilde{\rho}^k)$ is a quasi-statistical structure if and only if $\gamma\neq0,~~~\alpha=\beta=0$.
\end{thm}
For $G_3$, we have
\begin{align}
T^c(\widetilde{e}_1,\widetilde{e}_2)=\gamma \widetilde{e}_3,~~~T^c(\widetilde{e}_1,\widetilde{e}_3)=(\beta-m_3)\widetilde{e}_2,~~~T^c(\widetilde{e}_2,\widetilde{e}_3)=(m_3-\alpha)\widetilde{e}_1.
\end{align}
\begin{align}
&\widetilde{\rho}^c(T^c(\widetilde{e}_1,\widetilde{e}_2),\widetilde{e}_1)=0,~~~\widetilde{\rho}^c(T^c(\widetilde{e}_1,\widetilde{e}_2),\widetilde{e}_2)=\widetilde{\rho}^c(T^c(\widetilde{e}_1,\widetilde{e}_2),\widetilde{e}_3)=0,\widetilde{\rho}^c(T^c(\widetilde{e}_1,\widetilde{e}_3),\widetilde{e}_1)=0,\\\nonumber
&\widetilde{\rho}^c(T^c(\widetilde{e}_1,\widetilde{e}_3),\widetilde{e}_2)=\gamma m_3(m_3-\beta),~~~\widetilde{\rho}^c(T^c(\widetilde{e}_1,\widetilde{e}_3),\widetilde{e}_1)=\gamma m_3(\alpha-m_3),\\\nonumber
&\widetilde{\rho}^c(T^c(\widetilde{e}_2,\widetilde{e}_3),\widetilde{e}_1)=\widetilde{\rho}^c(T^c(\widetilde{e}_2,\widetilde{e}_3),\widetilde{e}_1)=\widetilde{\rho}^c(T^c(\widetilde{e}_2,\widetilde{e}_3),\widetilde{e}_3)=0.
\end{align}
Then, if $(G_3,\nabla^c,\widetilde{\rho}^c)$ is a quasi-statistical structure, by (3.2) and (3.3), we have the following two equations:
\begin{eqnarray}
        \begin{cases}
       m_3\gamma(m_3-\beta)=0 \\[2pt]
             m_3\gamma(\alpha-m_3)=0\\[2pt]
       \end{cases}
\end{eqnarray}
By solving (5.17), we get\\
\begin{thm}
 $(G_3,\nabla^c,\widetilde{\rho}^c)$ is a quasi-statistical structure if and only if
  \begin{eqnarray*}
  &&(1)\gamma=0,\\
  &&(2)\gamma\neq0,~~~\alpha+\beta-\gamma=0.
  \end{eqnarray*}
\end{thm}
Similarly,
\begin{align}
T^k(\widetilde{e}_1,\widetilde{e}_2)=\gamma \widetilde{e}_3,~~~T^k(\widetilde{e}_1,\widetilde{e}_3)=T^k(\widetilde{e}_2,\widetilde{e}_3)=0.
\end{align}
\begin{align}
&\widetilde{\rho}^k(T^k(\widetilde{e}_1,\widetilde{e}_2),\widetilde{e}_j)=\widetilde{\rho}^k(T^k(\widetilde{e}_1,\widetilde{e}_3),\widetilde{e}_j)=\widetilde{\rho}^k(T^k(\widetilde{e}_2,\widetilde{e}_3),\widetilde{e}_j)=0.
\end{align}
where $1\leq j\leq3$.\\
Then, we get
\begin{thm}
$(G_3,\nabla^k,\widetilde{\rho}^k)$ is a quasi-statistical structure.
\end{thm}
For $G_4$, we have
\begin{align}
T^c(\widetilde{e}_1,\widetilde{e}_2)=(\beta-2\eta) \widetilde{e}_3,~~~T^c(\widetilde{e}_1,\widetilde{e}_3)=(\beta-n_3)\widetilde{e}_2-\widetilde{e}_3,~~~T^c(\widetilde{e}_2,\widetilde{e}_3)=(n_3-\alpha)\widetilde{e}_1.
\end{align}
\begin{align}
&\widetilde{\rho}^c(T^c(\widetilde{e}_1,\widetilde{e}_2),\widetilde{e}_1)=0,~~~\widetilde{\rho}^c(T^c(\widetilde{e}_1,\widetilde{e}_2),\widetilde{e}_2)=\frac{(\beta-2\eta)(n_3-\beta)}{2},~~~\widetilde{\rho}^c(T^c(\widetilde{e}_1,\widetilde{e}_2),\widetilde{e}_3)=0,\nonumber\\
&\widetilde{\rho}^c(T^c(\widetilde{e}_1,\widetilde{e}_3),\widetilde{e}_1)=0,~~~\widetilde{\rho}^c(T^c(\widetilde{e}_1,\widetilde{e}_3),\widetilde{e}_2)=(\beta-n_3)[(2\eta-\beta)n_3-\frac{1}{2}],~~~\widetilde{\rho}^c(T^c(\widetilde{e}_1,\widetilde{e}_3),\widetilde{e}_1)=-\frac{(n_3-\beta)^2}{2},\\\nonumber
&\widetilde{\rho}^c(T^c(\widetilde{e}_2,\widetilde{e}_3),\widetilde{e}_1)=(n_3-\alpha)[(2\eta-\beta)n_3-1],~~~\widetilde{\rho}^c(T^c(\widetilde{e}_2,\widetilde{e}_3),\widetilde{e}_2)=\widetilde{\rho}^c(T^c(\widetilde{e}_2,\widetilde{e}_3),\widetilde{e}_3)=0.
\end{align}
Then, if $(G_4,\nabla^c,\widetilde{\rho}^c)$ is a quasi-statistical structure, by (3.2) and (3.3), we have the following five equations:
\begin{eqnarray}
        \begin{cases}
       \frac{(\beta-2\eta)(n_3-\beta)}{2}=0 \\[2pt]
       \frac{n_3-\beta}{2}=0 \\[2pt]
       (\beta-n_3)[(2\eta-\beta)n_3-\frac{1}{2}]=0 \\[2pt]
       \frac{(2n_3-\beta)(n_3-\beta)}{2}=0 \\[2pt]
             \frac{\beta-n_3}{2}+(n_3-\alpha)[(2\eta-\beta)n_3-1]=0\\[2pt]
       \end{cases}
\end{eqnarray}
By solving (5.22), we get\\
\begin{thm}
 $(G_4,\nabla^c,\widetilde{\rho}^c)$ is a quasi-statistical structure if and only if
  \begin{eqnarray*}
  &&(1)\alpha=\beta=2\eta,\\
  &&(2)\alpha=0,~~~\beta=\eta.
  \end{eqnarray*}
\end{thm}
Similarly,
\begin{align}
T^k(\widetilde{e}_1,\widetilde{e}_2)=(\beta-2\eta) \widetilde{e}_3,~~~T^k(\widetilde{e}_1,\widetilde{e}_3)=T^k(\widetilde{e}_2,\widetilde{e}_3)=0.
\end{align}
\begin{align}
&\widetilde{\rho}^k(T^k(\widetilde{e}_1,\widetilde{e}_2),\widetilde{e}_1)=0,~~~\widetilde{\rho}^k(T^k(\widetilde{e}_1,\widetilde{e}_2),\widetilde{e}_2)=\frac{\alpha(\beta-2\eta)}{2},~~~\widetilde{\rho}^k(T^k(\widetilde{e}_1,\widetilde{e}_2),\widetilde{e}_3)=0,\\\nonumber
&\widetilde{\rho}^k(T^k(\widetilde{e}_1,\widetilde{e}_3),\widetilde{e}_j)=0,~~~\widetilde{\rho}^k(T^k(\widetilde{e}_2,\widetilde{e}_3),\widetilde{e}_j)=0.
\end{align}
where $1\leq j\leq3$.\\
Then, if $(G_4,\nabla^k,\widetilde{\rho}^k)$ is a quasi-statistical structure, by (3.2) and (3.3), we have the following three equations:
\begin{eqnarray}
        \begin{cases}
       (2\eta-\beta)(\frac{\alpha}{2}-\beta)=0 \\[2pt]
              \frac{\alpha}{2}-\beta=0 \\[2pt]
      \frac{\alpha\beta}{2}=0\\[2pt]
       \end{cases}
\end{eqnarray}
By solving (5.25), we get\\
\begin{thm}
 $(G_4,\nabla^k,\widetilde{\rho}^k)$ is a quasi-statistical structure if and only if $\alpha=\beta=0$.
\end{thm}
For $G_5$, we have
\begin{align}
T^c(\widetilde{e}_1,\widetilde{e}_2)=0,~~~T^c(\widetilde{e}_1,\widetilde{e}_3)=-\alpha\widetilde{e}_1-\frac{\beta+\gamma}{2}\widetilde{e}_2,~~~T^c(\widetilde{e}_2,\widetilde{e}_3)=-\frac{\beta+\gamma}{2}\widetilde{e}_1-\delta\widetilde{e}_2.
\end{align}
\begin{align}
&\widetilde{\rho}^c(T^c(\widetilde{e}_1,\widetilde{e}_2),\widetilde{e}_j)=\widetilde{\rho}^c(T^c(\widetilde{e}_1,\widetilde{e}_3),\widetilde{e}_j)=\widetilde{\rho}^c(T^c(\widetilde{e}_2,\widetilde{e}_3),\widetilde{e}_j)=0.
\end{align}
where $1 \leq j\leq 3$.\\
Then, we get
\begin{thm}
$(G_5,\nabla^c,\widetilde{\rho}^c)$ is a quasi-statistical structure.
\end{thm}
Similarly,
\begin{align}
T^k(\widetilde{e}_1,\widetilde{e}_2)=T^k(\widetilde{e}_1,\widetilde{e}_3)=T^k(\widetilde{e}_2,\widetilde{e}_3)=0.
\end{align}
\begin{align}
&\widetilde{\rho}^k(T^k(\widetilde{e}_1,\widetilde{e}_2),\widetilde{e}_j)=\widetilde{\rho}^k(T^k(\widetilde{e}_1,\widetilde{e}_3),\widetilde{e}_j)=\widetilde{\rho}^k(T^k(\widetilde{e}_2,\widetilde{e}_3),\widetilde{e}_j)=0.
\end{align}
where $1\leq j\leq3$.\\
Then, we get
\begin{thm}
$(G_5,\nabla^k,\widetilde{\rho}^k)$ is a quasi-statistical structure.
\end{thm}
For $G_6$, we have
\begin{align}
T^c(\widetilde{e}_1,\widetilde{e}_2)=-\beta \widetilde{e}_3,~~~T^c(\widetilde{e}_1,\widetilde{e}_3)=-\frac{\beta+\gamma}{2}\widetilde{e}_2-\delta \widetilde{e}_3,~~~T^c(\widetilde{e}_2,\widetilde{e}_3)=\frac{\beta-\gamma}{2}\widetilde{e}_1.
\end{align}
\begin{align}
&\widetilde{\rho}^c(T^c(\widetilde{e}_1,\widetilde{e}_2),\widetilde{e}_1)=0,~~~\widetilde{\rho}^c(T^c(\widetilde{e}_1,\widetilde{e}_2),\widetilde{e}_2)=\frac{\beta}{2}[\alpha\gamma-\frac{\delta}{2}(\beta-\gamma)],~~~\widetilde{\rho}^c(T^c(\widetilde{e}_1,\widetilde{e}_2),\widetilde{e}_3)=0,\nonumber\\
&\widetilde{\rho}^c(T^c(\widetilde{e}_1,\widetilde{e}_3),\widetilde{e}_1)=0,~~~\widetilde{\rho}^c(T^c(\widetilde{e}_1,\widetilde{e}_3),\widetilde{e}_2)=\frac{\beta+\gamma}{2}[\alpha^2-\frac{1}{2}\beta(\beta-\gamma)]+\frac{\delta}{2}[\alpha\gamma-\frac{\delta}{2}(\beta-\gamma)],\nonumber\\
&\widetilde{\rho}^c(T^c(\widetilde{e}_1,\widetilde{e}_3),\widetilde{e}_3)=\frac{\beta+\gamma}{4}[\alpha\gamma-\frac{\delta}{2}(\beta-\gamma)],~~~\widetilde{\rho}^c(T^c(\widetilde{e}_2,\widetilde{e}_3),\widetilde{e}_1)=\frac{\beta-\gamma}{2}[\frac{\beta}{2}(\beta-\gamma)-\alpha^2],\nonumber\\
&\widetilde{\rho}^c(T^c(\widetilde{e}_2,\widetilde{e}_3),\widetilde{e}_1)=\widetilde{\rho}^c(T^c(\widetilde{e}_2,\widetilde{e}_3),\widetilde{e}_3)=0.
\end{align}
Then, if $(G_6,\nabla^c,\widetilde{\rho}^c)$ is a quasi-statistical structure, by (3.2) and (3.3), we have the following five equations:
\begin{eqnarray}
        \begin{cases}
      \frac{\beta}{2}[\alpha\gamma-\frac{1}{2}\delta(\beta-\gamma)]=0 \\[2pt]
        \frac{\alpha}{2}[\alpha\gamma-\frac{1}{2}\delta(\beta-\gamma)]=0\\[2pt]
              \frac{\beta+\gamma}{2}[\alpha^2-\frac{1}{2}\beta(\beta-\gamma)]+\frac{\delta}{2}[\alpha\gamma-\frac{1}{2}\delta(\beta-\gamma)]=0\\[2pt]
 \frac{\gamma}{2}[\alpha\gamma-\frac{1}{2}\delta(\beta-\gamma)]=0\\[2pt]        \frac{\alpha}{2}[\alpha\gamma-\frac{1}{2}\delta(\beta-\gamma)]+\frac{\beta-\gamma}{2}[\alpha^2-\frac{1}{2}\beta(\beta-\gamma)]=0\\[2pt]
       \end{cases}
\end{eqnarray}
By solving (5.32), we get\\
\begin{thm}
 $(G_6,\nabla^c,\widetilde{\rho}^c)$ is a quasi-statistical structure if and only if
 \begin{eqnarray*}
 &&(1)\alpha=\beta=\gamma=0,~~~\delta\neq0,\\
 &&(2)\alpha\neq0,~~~\delta=\gamma=0,~~~2\alpha^2=\beta^2,\\
 &&(3)\alpha\neq0,~~~\beta=\gamma=0.
 \end{eqnarray*}
\end{thm}
Similarly,
\begin{align}
T^k(\widetilde{e}_1,\widetilde{e}_2)=-\beta \widetilde{e}_3,~~~T^k(\widetilde{e}_1,\widetilde{e}_3)=T^k(\widetilde{e}_2,\widetilde{e}_3)=0.
\end{align}
\begin{align}
&\widetilde{\rho}^k(T^k(\widetilde{e}_1,\widetilde{e}_2),\widetilde{e}_j)=\widetilde{\rho}^k(T^k(\widetilde{e}_1,\widetilde{e}_3),\widetilde{e}_j)=\widetilde{\rho}^k(T^k(\widetilde{e}_2,\widetilde{e}_3),\widetilde{e}_j)=0.
\end{align}
where $1\leq j\leq3$.\\
Then, if $(G_6,\nabla^k,\widetilde{\rho}^k)$ is a quasi-statistical structure, by (3.2) and (3.3), we have the following two equations:
\begin{eqnarray}
        \begin{cases}
       \alpha\beta\gamma=0 \\[2pt]
             \alpha^2\gamma=0\\[2pt]
       \end{cases}
\end{eqnarray}
By solving (5.35), we get\\
\begin{thm}
 $(G_6,\nabla^k,\widetilde{\rho}^k)$ is a quasi-statistical structure if and only if
 \begin{eqnarray*}
 &&(1)\alpha=\beta=0,~~~\delta\neq0,\\
 &&(2)\alpha\neq0,~~~\beta\delta=\gamma=0.
 \end{eqnarray*}
\end{thm}
For $G_7$, we have
\begin{align}
T^c(\widetilde{e}_1,\widetilde{e}_2)=\beta \widetilde{e}_3,~~~T^c(\widetilde{e}_1,\widetilde{e}_3)=-\alpha\widetilde{e}_1-\frac{\gamma}{2}\widetilde{e}_2-\beta \widetilde{e}_3,~~~T^c(\widetilde{e}_2,\widetilde{e}_3)=-(\beta+\frac{\gamma}{2})\widetilde{e}_1-\delta\widetilde{e}_2-\delta\widetilde{e}_3.
\end{align}
\begin{align}
&\widetilde{\rho}^c(T^c(\widetilde{e}_1,\widetilde{e}_2),\widetilde{e}_1)=-\frac{\beta}{2}(\alpha\gamma+\frac{\delta\gamma}{2}),~~~\widetilde{\rho}^c(T^c(\widetilde{e}_1,\widetilde{e}_2),\widetilde{e}_2)=\frac{\beta}{2}(\alpha^2+\frac{\beta\gamma}{2}),~~~\widetilde{\rho}^B(T^B(\widetilde{e}_1,\widetilde{e}_2),\widetilde{e}_3)=0,\\\nonumber
&\widetilde{\rho}^c(T^c(\widetilde{e}_1,\widetilde{e}_3),\widetilde{e}_1)=\alpha^3+\alpha\beta\gamma+\frac{\beta\delta\gamma}{4},~~~\widetilde{\rho}^c(T^c(\widetilde{e}_1,\widetilde{e}_3),\widetilde{e}_2)=(\frac{\alpha^2}{2}+\frac{\beta\gamma}{4})(\gamma-\beta),\\\nonumber
&\widetilde{\rho}^c(T^c(\widetilde{e}_1,\widetilde{e}_3),\widetilde{e}_1)=\frac{\alpha\gamma(\alpha+\delta)}{4}-\frac{\beta\gamma^2}{8}),~~~\widetilde{\rho}^c(T^c(\widetilde{e}_2,\widetilde{e}_3),\widetilde{e}_1)=(\beta+\frac{\gamma}{2})(\alpha^2+\frac{\beta\gamma}{2})+\frac{\delta\gamma}{2}(\alpha+\frac{\delta}{2}),\\\nonumber
&\widetilde{\rho}^c(T^c(\widetilde{e}_2,\widetilde{e}_3),\widetilde{e}_1)=\frac{\delta}{2}(\alpha^2+\frac{\beta\gamma}{2}),~~~\widetilde{\rho}^c(T^c(\widetilde{e}_2,\widetilde{e}_3),\widetilde{e}_3)=\frac{\alpha\beta\gamma-\alpha^2\delta}{2}+\frac{\alpha\gamma^2}{4}+\frac{\delta\gamma^2}{8}.
\end{align}
Then, if $(G_7,\nabla^c,\widetilde{\rho}^c)$ is a quasi-statistical structure, by (3.2) and (3.3), we have the following nine equations:
\begin{eqnarray}
        \begin{cases}
       \frac{\beta\gamma}{2}(\alpha+\frac{\delta}{2})=0 \\[2pt]
       \frac{\beta}{2}(\alpha^2+\frac{\beta\gamma}{2})=0 \\[2pt]
       \frac{\beta}{2}(\alpha^2+\frac{\beta\gamma}{2})-\frac{\alpha\gamma}{2}(\alpha+\frac{\delta}{2})=0  \\[2pt]
       \frac{\alpha^3}{2}+\frac{3\alpha\beta\gamma+\beta\delta\gamma}{4})=0 \\[2pt]
      (\gamma-\beta)(\frac{\alpha^2+\frac{\beta\gamma}{4}})-\frac{\alpha\gamma}{2}(\alpha+\frac{\delta}{2})=0 \\[2pt]
       \frac{1}{2}(\frac{\gamma}{2}-\beta)(\alpha^2+\frac{\beta\gamma}{2})+\frac{\alpha\delta\gamma+\alpha^2\gamma}{4})-\frac{\beta\gamma^2}{8}=0\\[2pt]
       \frac{\delta\gamma}{2}(\alpha+\frac{\delta}{2}\frac{\gamma-\beta}{2}(\alpha^2+\frac{\beta\gamma}{2})=0 \\[2pt]
      \frac{\delta}{2}(\alpha^2+\frac{\beta\gamma}{2})-\frac{\beta\gamma}{2}(\alpha+\frac{\delta}{2})=0\\[2pt]
        \frac{\alpha\beta\gamma-\alpha^2\delta}{2}+\frac{\alpha\gamma^2}{4}+\frac{\delta\gamma^2}{8}-\frac{1}{2}(\beta-\frac{\gamma}{2})(\alpha\gamma+\frac{\delta\gamma}{2})=0\\[2pt]
       \end{cases}
\end{eqnarray}
By solving (5.38), we get\\
\begin{thm}
 $(G_7,\nabla^c,\widetilde{\rho}^c)$ is a quasi-statistical structure if and only if $\delta\neq0,~~~\alpha=\gamma=0$.
\end{thm}
Similarly,
\begin{align}
T^k(\widetilde{e}_1,\widetilde{e}_2)=\beta \widetilde{e}_3,~~~T^k(\widetilde{e}_1,\widetilde{e}_3)=T^k(\widetilde{e}_2,\widetilde{e}_3)=0.
\end{align}
\begin{align}
&\widetilde{\rho}^k(T^k(\widetilde{e}_1,\widetilde{e}_2),\widetilde{e}_1)=\beta^2(\alpha+\delta),~~~\widetilde{\rho}^k(T^k(\widetilde{e}_1,\widetilde{e}_2),\widetilde{e}_2)=\frac{\beta}{2}(\alpha\delta+\beta\gamma+2\delta^2),~~~\widetilde{\rho}^k(T^k(\widetilde{e}_1,\widetilde{e}_2),\widetilde{e}_3)=0,\\\nonumber
&\widetilde{\rho}^k(T^k(\widetilde{e}_1,\widetilde{e}_3),\widetilde{e}_j)=0,~~~\widetilde{\rho}^k(T^k(\widetilde{e}_2,\widetilde{e}_3),\widetilde{e}_j)=0.
\end{align}
where $1\leq j\leq3$.\\
Then, if $(G_7,\nabla^k,\widetilde{\rho}^k)$ is a quasi-statistical structure, by (3.2) and (3.3), we have the following nine equations:
\begin{eqnarray}
        \begin{cases}
       \beta(\alpha\gamma+\beta\delta+2\beta\delta)=0 \\[2pt]
              \beta(\delta^2-\alpha^2-\beta^2+\frac{3\alpha\delta-\beta\gamma}{2})=0 \\[2pt]
       \beta(\alpha+\delta)^2=0 \\[2pt]
       2\alpha^3-\alpha\delta^2-2\beta^2\delta-\frac{\alpha\beta\gamma+\alpha^2\delta}{2}=0 \\[2pt]
              \beta(\alpha^2+\beta^2+\frac{\beta\gamma+\delta^2}{2})+\alpha^2\gamma=0 \\[2pt]
      -\alpha\beta(\alpha+\delta)-\frac{\beta}{2}(\beta\gamma+\alpha\delta+2\delta^2)=0\\[2pt]
             \alpha^2\gamma+\beta^3+\frac{3\alpha^2\beta+\beta^2\gamma-3\alpha\beta\delta-5\beta\delta^2}{2}=0 \\[2pt]
       \alpha\beta\gamma-\delta^3+3\beta^2\delta+2\alpha^2\delta+\frac{\beta\delta\gamma-\alpha\delta^2}{2}=0 \\[2pt]
       \alpha\beta\gamma+\delta^3+\frac{3\beta\delta\gamma+\alpha\delta^2}{2}=0 \\[2pt]
       \end{cases}
\end{eqnarray}
By solving (5.41), we get $\alpha=\delta=0$, there is a contradiction. So\\
\begin{thm}
 $(G_7,\nabla^k,\widetilde{\rho}^k)$ is not a quasi-statistical structure.
\end{thm}
\section*{Acknowledgements}
The second author was supported in part by  NSFC No.11771070. The authors are deeply grateful to the referees for their valuable comments and helpful suggestions.

\end{document}